\title{Descendant invariants and characteristic numbers}
\author{Tom Graber
\and Joachim Kock
\and Rahul Pandharipande}
\date{}

\documentclass[12pt]{article}
\usepackage{amsmath}
\usepackage{amssymb}
\usepackage{latexsym}
\input{diagrams}
\newcommand{\hovedfont}{\normalfont\bfseries}
\usepackage{theorem}
\theorembodyfont{\normalfont\itshape}
\theoremheaderfont{\hovedfont}
	\theoremstyle{change}
\newtheorem{lemma}{Lemma.}[subsection]
\newtheorem{satz}[lemma]{Theorem.}
\newtheorem{prop}[lemma]{Proposition.}
\newtheorem{cor}[lemma]{Corollary.}
	\theorembodyfont{\normalfont}
\newtheorem{eks}[lemma]{Example.}

\newtheorem{BM}[lemma]{Remark.}
\newtheorem{taller}[lemma]{$\!\!$}
	\newenvironment{blanko}[1]%
{\begin{taller}{\hovedfont #1}\normalfont}%
{\end{taller}}

	\newenvironment{dem}%
{\begin{list}{\em Proof. }%
{\setlength{\labelsep}{0mm}\setlength{\leftmargin}{0mm}%
\setlength{\labelwidth}{0mm}\setlength{\listparindent}{\parindent}%
\setlength{\parsep}{\parskip}\setlength{\partopsep}{0mm}}%
\item}{\qed\end{list}}
	\newenvironment{dem*}[1]%
{\begin{list}{\em #1 }%
{\setlength{\labelsep}{0mm}\setlength{\leftmargin}{0mm}%
\setlength{\labelwidth}{0mm}\setlength{\listparindent}{\parindent}%
\setlength{\parsep}{\parskip}\setlength{\partopsep}{0mm}}%
\item}{\qed\end{list}}
	\newenvironment{blanko*}[1]%
{\begin{list}{\bf {#1} }%
{\setlength{\labelsep}{0mm}\setlength{\leftmargin}{0mm}%
\setlength{\labelwidth}{0mm}\setlength{\listparindent}{\parindent}%
\setlength{\parsep}{\parskip}\setlength{\partopsep}{0mm}}%
\item}{\end{list}}
\newcounter{dummycounter}
\newenvironment{punkt-i}%
{%
	\begin{list}%
	{(\roman{dummycounter})}%
	{\usecounter{dummycounter}%
	\setlength{\itemsep}{0em}\setlength{\parsep}{0em}\setlength{\topsep}{0em}%
	\setlength{\itemindent}{0em}\setlength{\labelwidth}{1.8em}%
	\setlength{\labelsep}{0.6em}\setlength{\leftmargin}{2.4em}}%
}%
{\end{list}}
\renewcommand{\epsilon}{\varepsilon}
\makeatletter
\renewcommand{\ldots}{\relax\ifmmode\ldotp\ldotp\ldotp%
\else$\m@th\ldotp\ldotp\ldotp\ $\fi}
\makeatother
\newcommand{\brac}[1]{\,\langle \, {#1} \, \rangle\,}
\newcommand{\bbrac}[1]%
{\,\langle\!\langle \, {#1} \, \rangle\!\rangle\,}
\newcommand{\Brac}[1]{\,\langle \ {#1} \ \rangle\,}
\newcommand{\Bbrac}[1]%
{\,\langle\!\langle \ {#1} \ \rangle\!\rangle\,}
\newcommand{\cp}{ {\scriptstyle \;\cup\;} }
\newcommand{\til}{\tilde}

\newcommand{\smallbinom}[2]%
{{\textstyle \binom{#1}{#2}}}
\newcommand{\smallprod}[2]%
{\overset{#2}{\underset{#1}{\textstyle{\prod}}}} 
\newcommand{\smallsum}[2]%
{\overset{#2}{\underset{#1}{\textstyle{\sum}}}} 
\newcommand{\bigprod}[2]%
{\overset{#2}{\underset{#1}{\text{{\Large\(\prod\)}}}}}
\newcommand{\bigsum}[2]%
{\overset{#2}{\underset{#1}{\text{{\Large\(\sum\)}}}}}

\newcommand{\bigfatgreek}[1]{\boldsymbol{#1}} 
 
\newcommand{\psiclass}{\bigfatgreek\psi} 
\newcommand{\mpsiclass}{\ov\psiclass{}} 
\newcommand{\deltaclass}{\bigfatgreek\delta} 
\newcommand{\boundary}{\bigfatgreek\xi}

\newcommand{\XX}{\mathfrak{X}}
\providecommand{\qed}%
{\hspace*{\fill}\nolinebreak[1]\hspace*{\fill} $\Box$}

\renewcommand{\d}{\partial}

\providecommand{\norm}[1]{\left| {#1}\right|}

\newcommand{\pil}{\rightarrow}
\newcommand{\langpil}{\longrightarrow}

\newcommand{\into}{\hookrightarrow}
\newcommand{\isopil}%
{\stackrel{\raisebox{0.1ex}[0ex][0ex]{\(\sim\)}}%
{\raisebox{-0.15ex}[0.28ex]{\(\pil\)}}}

\newcommand{\shortsetminus}%
{\,\raisebox{1pt}{\ensuremath{\mathbb{r}}\,}}
\newcommand{\kan}{\mbox{\Large \(\omega\)}}
\newcommand{\OO}{\ensuremath{{\mathcal O}}}

\newcommand{\upperstar}{^{\raisebox{-0.25ex}[0ex][0ex]{\(\ast\)}}}
\newcommand{\lowerstar}{_{\raisebox{-0.33ex}[-0.5ex][0ex]{\(\ast\)}}}
\newcommand{\df}%
{\: {\raisebox{0.255ex}{\normalfont\scriptsize :\!\!}}=}
\newcommand{\tensor}{\otimes}
\renewcommand{\dim}{\operatorname{dim}}

\newcommand{\Gr}{\operatorname{Gr}}

\newcommand{\brok}[2]{{\textstyle{\frac{#1}{#2}}}}

\newcommand{\ov}{\overline}

\newcommand{\fat}[1]{\mathbf{{#1}}}

\renewcommand{\P}	{\mathbb{P}}

\newcommand{\Z}		{\mathbb{Z}}
\newcommand{\Q}		{\mathbb{Q}}

\newcommand{\mtau}{\ov\tau}
\newcommand{\mbtau}{\overline{\boldsymbol{\tau}}}
\newcommand{\vdim}{\operatorname{vdim}}

\newcommand{\M}{{\overline M}}
\newcommand{\dual}{^{\vee}}          
\newcommand{\virtual}{^{\text{vir}}} 
\newcommand{\evalmap}{\nu}              
\newcommand{\LL}{\mathbb{L}}         
\newcommand{\Hodge}{\mathbb{E}}      
\newcommand{\EE}{{\mathcal{E}}}      

\newcommand{\incidencelocus}[1]{I_{{#1}}}
\newcommand{\pointedincidencelocus}[2]{\M{}_{\{{#1}\}}^{{#2}}}
\newcommand{\tangencylocus}{T}
\newcommand{\pointedtangencylocus}{Z}

\newcommand{\incidenceclass}{\mathcal{I}} 
\newcommand{\tangencyclass}{\mathcal{T}} 
\newcommand{\flagclass}{\mathcal{F}}

\newcommand{\rarr}{\rightarrow}

\newcommand{\sumo}{\oplus}

\newcommand{\soft}{{\em{easy }}}

\begin{document}
\maketitle 
\begin{abstract}
	On a stack of stable maps, the psi classes are modified by subtracting
	certain boundary divisors.  These modified psi classes are compatible with
	forgetful morphisms, and are well-suited to enumerative geometry: tangency
	conditions allow simple expressions in terms of modified psi classes. 
	Topological recursion relations are established among their top products in
	genus zero, yielding effective recursions for characteristic numbers of
	rational curves in any projective homogeneous variety.  In higher genus, the
	obtained numbers are only virtual, due to contributions from spurious
	components of the space of maps.  For the projective plane, the necessary
	corrections are determined in genus $1$ and $2$ to give the characteristic
	numbers in these cases.
\end{abstract}

\section*{Introduction}

\begin{blanko*}{Gromov-Witten invariants and enumerative geometry.}
	\ \sloppy The Gromov-Witten in\-variants of a nonsingular algebraic variety
	$X$ are defined via integrals over associated moduli stacks of maps from
	curves to $X$ (cf.~Kontsevich-Manin~\cite{KM:9402},
	Ruan-Tian~\cite{Ruan-Tian}), Behrend~\cite{Behrend}).  One of the motivating
	properties of Gromov-Witten invariants is their connection to enumerative
	geometry.  The simplest relationship occurs when one considers the genus $0$
	Gromov-Witten invariants of a compact homogeneous variety.  These invariants
	are exactly the solutions to associated problems in classical enumerative
	geometry.
	\fussy

Even more is true in case the target is $\P^2$.  It is not hard to prove that
the positive degree Gromov-Witten invariants of $\P^2$ are precisely the
Severi degrees: the numbers of degree $d$, genus $g$ plane curves passing
through $3d+g-1$ general points.  More generally, one can ask for the number of
degree $d$, genus $g$ plane curves passing through $a$ general points and
tangent to $b$ general lines ($a+b= 3d+g-1$).  These are the {\em characteristic
numbers} of plane curves.  These numbers were of great interest to the classical
algebraic geometers --- they were tabulated in the nineteenth century for degrees
up to 4 by H.~Zeuthen~\cite{Zeuthen} and H.~Schubert~\cite{Schubert}.  Zeuthen's
degeneration methods have been recently reinterpreted and vindicated in the
context of stable maps by R.~Vakil in \cite{Vakil:9812}.

A new method of attacking the characteristic number problem is pursued here, 
leading to concise solutions to the problem of characteristic numbers of 
rational curves in any homogeneous variety, and for plane curves of genus $0$, 
$1$, and $2$.
The idea is to interpret the tangency conditions in terms of gravitational
descendant integrals over the moduli space of maps.  
\end{blanko*}

\begin{blanko*}{Gravitational descendants.}
	An important generalization of the Gromov-Witten invariants is to include the
	tautological cotangent line classes (psi classes) in the integrals, which in
	theoretical physics corresponds to introducing gravity into the topological
	sigma model, cf.~Witten~\cite{Witten}.  These {\em gravitational descendants}
	are now a central notion in Gromov-Witten theory: they appear as fundamental
	solutions to Givental's quantum differential equation (see for example
	\cite{Pand:9806}), and they are the subject of the 
Virasoro conjecture (see
	Getzler~\cite{Getzler:9812} for a survey.)
%
\end{blanko*}

\begin{blanko*}{Modified psi classes and descendants.}	
	This work presents a connection between gravitational quantum
	coho\-mology and enumerative geometry.  The gravitational descendants are not
	in general enumerative, but  a modification of the psi class is introduced,
	which is well-suited to enumerative geometry, and to tangency conditions in
	particular.  For example, if $X$ is a compact homogeneous variety, the
	characteristic numbers of rational curves in $X$ are top products of modified
	psi classes and evaluation classes.
	
	The first two sections are devoted to foundational issues concerning modified
	psi classes, as well as their companion classes (called
	diagonal classes).  Unlike the tautological psi classes, the modified ones
	pull-back to the boundary in a non-trivial way, giving rise to diagonal
	classes as correction terms.  Topological recursion relations (TRR) do exist
	in genus $0$ and $1$ (as described in Section~\ref{Sec:toprec}), but their
	combinatorics is much messier than in the tautological case, due to the
	appearance of the diagonal classes.  For those invariants whose exponents on
	each modified psi class are at most one (called first enumerative
	descendants), which are the ones that arise in the characteristic number
	problem, the notion of the deformed metric of \cite{Kock:0006} clarifies the
	recursions and yields concise equations.

\end{blanko*}

\begin{blanko*}{Characteristic numbers of rational curves in a homogeneous variety.}
	The condition of being tangent to a hypersurface admits a simple expression
	in terms of modified psi classes.  In genus $0$, it is shown that top
	products of these expressions are exactly the characteristic numbers.  The
	characteristic number problem is then given a concise solution in terms of
	partial differential equations for their generating functions; these
	equations are simple coordinate changes of the TRR for the first enumerative
	descendants.
  
	Historically, the case $X=\P^2$ has been of special interest.  
A very pleasant
	solution to the characteristic number problem for rational plane curves
	follows as a special case of these equations.  Consider the full rational
	characteristic number potential,
	\begin{equation*}
	G^{(0)}= \sum_{d\geq 1} \exp(ds)\sum_{a,b,c}
	\frac{u^a}{a!}\frac{v^b}{b!}\frac{w^c}{c!}
	N^{(0)}_d\!(a,b,c) , 
	\end{equation*}
	where $N^{(0)}_d\!(a,b,c)$ is the number of degree $d$ rational plane curves
	passing through $a$ general points, tangent to $b$ general lines, and tangent
	to $c$ general lines at a specified point at each line.  The potential 
	satisfies the differential equation
\begin{equation}
\label{dfe}
G^{(0)}_{v s} =  G^{(0)}_{u s} - G^{(0)}_{u}
+\brok{1}{2} G^{(0)}_{s s}G^{(0)}_{s s}
+ 2v G^{(0)}_{s s}G^{(0)}_{u s}
+ (v^2+w) G^{(0)}_{u s}G^{(0)}_{u s}.
\end{equation}
	This equation expresses each characteristic number in terms of
	those with strictly fewer tangency conditions.  There is a similar equation
	for decrementing the number of flag conditions.  This and the
	WDVV-equations determine $G^{(0)}$ from the initial condition
	$N^{(0)}_1\!(2,0,0)=1$.

	In the case of rational plane curves, a partial computation of the
	characteristic numbers was given in Di
	Francesco-Itzykson~\cite{diFrancesco-Itzykson}.  A complete solution was
	given in Pandharipande~\cite{Pand:9504} via intersection theory on
	$\M_{0,n}(\P^2,d)$.  A different solution was obtained in
	Ernstr\"om-Kennedy~\cite{EK1} via an investigation of contact spaces.  These
	solutions are complicated by auxiliary elements in the recursions. 
	R.~Vakil~\cite{Vakil:9803} has shown (\ref{dfe}) may be derived from
	degeneration arguments together with formulas in \cite{Pand:9504}.
\end{blanko*}

\begin{blanko*}{Characteristic numbers of plane curves.}
	For the projective plane, even the higher genus Gromov-Witten invariants are
	enumerative: although the moduli spaces have components of excessive
	dimension, corresponding to unwanted degenerate maps, it turns out that the
	incidence conditions are too strong for the ill-behaved components to
	contribute.

	Unfortunately, when it comes to tangency conditions, this is no longer so,
	and the solutions obtained by applying the above techniques directly will in
	general be virtual numbers.  To arrive at the enumerative solutions,
	extraneous contributions must be removed.  These excess contributions vanish
	in genus $0$ and may be evaluated explicitly in genus $1$ and $2$:  In genus
	$1$, the virtual number receives a contribution from maps which are a rational
	curve with a contracted elliptic tail attached.  In genus $2$, the
	extraneous contributions come from elliptic curves with a contracted
	elliptic tail attached, rational curves with two contracted elliptic tails
	attached and curves which are the union of rational curve and a genus $2$
	double cover of a line.  For each of these types of curves, the corresponding
	potential is identified.  
	A solution of the characteristic
	number problem for plane curves is thus obtained in genus $0$, $1$, and $2$
	(for all degrees).
	
	In genus $1$, a differential equation can be derived for the characteristic
	number potential by applying the corrections to TRR. With notation similar
	to that of genus $0$, one equation reads
\begin{multline*}
G^{(1)}_v = G^{(1)}_u 
+  G^{(0)}_{ss} \cdot  G^{(1)}_s + 2v G^{(0)}_{us} \cdot  G^{(1)}_s 
+ 2v G^{(0)}_{ss} \cdot  G^{(1)}_u + (2v^2+2w) G^{(0)}_{us} \cdot  G^{(1)}_u 
 \\ +
\frac{1}{24}\bigg( G^{(0)}_{s s s} -3 G^{(0)}_{s s} + 2 G^{(0)}_s \
+ \ 2v( 2G^{(0)}_{uss} - 2G^{(0)}_{us} - G^{(0)}_{vs}) \
+ \ (2v^2 + 2w) (G^{(0)}_{uus} - G^{(0)}_{ws})\bigg) .
\end{multline*}
This equation has also been derived by R.~Vakil~\cite{Vakil:9803}
by degeneration methods and formulas from \cite{Pand:9504}.
Similarly there is an equation for decrementing the number of flag 
conditions.  

In genus $2$, it seems too messy to extract a single equation for the
characteristic numbers.  The virtual numbers can be computed using localization
formulas \cite{Graber-Pandharipande:9708}, or via the descendant relation of
Belorousski-Pandharipande~\cite{Belorousski-Pandharipande} combined with the
genus $2$ TRR of Getzler~\cite{Getzler:9801}.  The correction terms involve only
the rational and elliptic potentials and are quickly determined by the equations
given above.

In the short Section~\ref{Sec:P1P1}, it is shown how the techniques also solve
the characteristic number problem for (rational and) elliptic curves in
$\P^1\times\P^1$.
\end{blanko*}

All the algorithms have been implemented in maple.  Code or numerical data is 
available upon request.

\begin{blanko*}{Acknowledgements.} 
	The first author was supported by a Sloan dissertation year fellowship and an
	NSF postdoctoral fellowship.  The second author was supported by the Natural
	Science Research Council of Denmark and the Nordic Research Academy (NorFA). 
	Part of this work was carried out while the second author was visiting
	the California Institute of Technology to which he is thankful for exquisite
	hospitality.
\end{blanko*}

\section{Modified psi classes and descendant invariants}

\subsection{Preliminaries}

\begin{blanko}{Set-up.}
	Throughout we work over the field of complex numbers.  Let $X$ be a
	nonsingular projective variety; let $T_0,\ldots,T_r$ denote the elements of a
	homogeneous additive basis for the cohomology space $H\upperstar (X,\Q)$; let
	$g_{ij}= \int_X T_i \cp T_j$ denote the entries of the intersection pairing
	matrix; and let $(g^{ij})$ be the inverse matrix.

	Let $\M_{g,S}(X,\beta)$ denote the moduli stack of Kontsevich stable maps to
	$X$ representing the class $\beta\in H_2^+(X,\Z)$, of genus $g$, and with
	marking set $S=\{p_1,\ldots,p_n\}$.  (When the names of the marks are not
	important the stack will also be denoted $\M_{g,n}(X,\beta)$.)  For each mark
	$p_i$, let $\evalmap_i : \M_{g,S}(X,\beta) \pil X$ denote the evaluation
	morphism, which sends the class of a map $[S \pil C\stackrel{\mu}{\pil}X]$ to
	$\mu(p_i)\in X$.  Pull-backs of cohomology classes in $X$ along these
	morphisms are called {\em evaluation classes}.
	
	Let $\pi_0 : \M_{g,S\cup\{p_0\}}(X,\beta) \pil \M_{g,S}(X,\beta)$ be the
	forgetful morphism which consists in forgetting the extra mark $p_0$ (and
	stabilizing, by contracting any component that becomes unstable in the
	absence of $p_0$).  The diagram
	\begin{diagram}[tight,h=6ex,w=3ex,shortfall=1.3ex]
\phantom{\XX}	& \M_{g,S\cup\{p_0\}}(X,\beta) &&& \rTo^{\evalmap_0} &&& X	\\
\uTo[shortfall=3.2ex,scriptlabels]<{\sigma_i \ }
                   &	\dTo>{\ \pi_0} &        &&&& 	 \\
\phantom{B}        & \M_{g,S}(X,\beta)
\end{diagram}
is the universal family of stable maps over $\M_{g,S}(X,\beta)$. 
(Cf.~\cite{Behrend-Manin}, 4.6.)

Here $\sigma_i$ is the section corresponding to the mark $p_i$.  The image of
$\sigma_i$ is the (closure of the) locus of maps whose source curve has two
components, one of which is rational, carries just the two marks $p_i$ and $p_0$,
and is contracted by $\mu$:
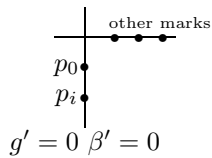
\begin{center}
\setlength{\unitlength}{4mm}
\footnotesize
\begin{picture}(4,5)(0,0)
   \put(1,1){\line(0,1){4}}
   \put(1,3){\circle*{.3}} \put(.4,3){\makebox(0,0){\(p_0\)}}
   \put(1,2){\circle*{.3}} \put(.4,2){\makebox(0,0){\(p_i\)}}
   \put(1,0.5){\makebox(0,0){\(g'=0\ \beta'=0 \)}}
   \put(0,4){\line(1,0){4}}
   \put(3.5,4.5){\makebox(0,0){\tiny other marks}}
	\put(2,4){\circle*{.25}}\put(2.8,4){\circle*{.25}}\put(3.6,4){\circle*{.25}}
\end{picture}
\end{center}
Let $D_{i,0}\in H^2(\M_{g,S\cup \{p_0\}}(X,\beta),\Q)$  denote the class of 
this Cartier divisor.
\end{blanko}

\begin{blanko}{Psi classes and descendants.}\label{psi}
	Let $\kan_{\pi_0}$ be the relative dualizing sheaf of $\pi_0$.  The $i$'th
	{\em cotangent line} of $\M_{g,S}(X,\beta)$ is the line bundle $\LL_i =
	\sigma_i\upperstar \kan_{\pi_0}$, whose fiber over a moduli point $[S \pil C
	\pil X]$ is the cotangent line $T_{p_i}\upperstar C$.  The $i$'th {\em psi
	class} is its first Chern class:
   $$
   \psiclass_i \df c_1(\LL_i) \in H^2(\M_{g,S}(X,\beta),\Q).
   $$
	
	The descendant invariants of $X$ are defined by integrals of products of
	evaluation classes and psi classes against the virtual fundamental class of
	$\M_{g,S}(X,\beta)$, (see Li-Tian~\cite{Li-Tian} and
	Behrend-Fantechi~\cite{Behrend-Fantechi}).  The following notation introduced
	by E.~Witten~\cite{Witten} has become standard.
$$
\brac{\tau_{a_1}(\gamma_1)\cdots \tau_{a_n}(\gamma_n)}\!_{g,\beta}^X \df
\int \psiclass_1^{a_1}\cp\evalmap_1\upperstar (\gamma_1) \cp\cdots\cp \psiclass_n^{a_n}
\cp\evalmap_n\upperstar (\gamma_n) 
\cap [\M_{g,n}(X,\beta) ]\virtual,
$$
where $\gamma_i \in H\upperstar (X,\Q)$ and the $a_i$'s are non-negative
integers.  The invariants are defined to vanish unless the dimension of the
integrand is equal to the expected dimension $(\dim X - 3)(1-g) + \int_\beta
c_1(T_X) + n$.  When the $a_i$'s are all zero, the descendants
specialize to the Gromov-Witten invariants of $X$.

For projective homogeneous varieties, the formulas of
virtual localization of \cite{Graber-Pandharipande:9708} determine the
descendants in any genus, for all degrees.  
\end{blanko}

\subsection{Modified psi classes and descendant invariants}

In enumerative geometry, the psi classes are not the most convenient classes;
their deficiency stems from the fact that they are not compatible with 
pull-back along forgetful morphisms. Indeed the formula reads
\begin{equation}\label{pull-psi}
   \pi_0\upperstar \psiclass_i = \psiclass_i - D_{i,0}.
\end{equation}
A modification of the psi class is introduced here, for enumerative purposes.

\begin{blanko}{Modified psi classes.}\label{mpsi}
	Suppose $\beta>0$ or $g>0$.
	For each mark $p_i \in S$, let 
	\begin{equation}\label{forgg}
	\hat\pi_i :\M_{g,S}(X,\beta) \pil \M_{g,\{p_i\}}(X,\beta)
	\end{equation}
	be the morphisms that forgets all marks but $p_i$.   The {\em modified
	psi} class on $\M_{g,S}(X,\beta)$ is by definition
	$$
   \mpsiclass_i \df \hat\pi_i\upperstar \psiclass_i.
   $$
	(The modified psi class is not defined for $(g,\beta)=(0,0)$; see however
	\cite{Kock:0006}, 4.4, for a discussion of this case.)  By construction, the
	modified psi classes are compatible with pull-backs along forgetful morphisms
	$\pi_0 : \M_{g,S\cup\{p_0\}}(X,\beta) \pil \M_{g,S}(X,\beta)$ (just as
	evaluation classes), in the sense that
   $$
   \pi_0 \upperstar \mpsiclass_i = \mpsiclass_i.
   $$
\end{blanko}

\begin{BM}
	The modified psi class defined here must not be confused with the phi class 
	studied in Kontsevich-Manin~\cite{KM:9708}: their class (which is also a 
	boundary modification of the psi class) is defined on $\M_{g,n}(X,\beta)$
	as the pull-back of the psi 
	class on $\M_{g,n}$ via the absolute stabilization morphism.
\end{BM}
\begin{blanko}{Enumerative descendants.}\label{edesc}
   We define an {\em enumerative descendant} to be a top product of modified 
   psi classes and evaluation classes, and agree on the 
   notation
   $$
   \Brac{\mtau_{a_1}(\gamma_1)\cdots \mtau_{a_n}(\gamma_n)}\!_{g,\beta}^X \df
\int \mpsiclass{}_1^{a_1}\cp\evalmap_1\upperstar (\gamma_1) \cp\cdots \cp
\mpsiclass{}_n^{a_n}\cp\evalmap_n\upperstar (\gamma_n)
\ \cap [\M_{g,n}(X,\beta) ]\virtual.
   $$
\end{blanko}

\begin{blanko}{String, dilaton, and divisor equations.}
	Since the modified psi classes as well as the evaluation classes are
	compatible with pull-back along forgetful morphisms, the projection formula
	readily implies the following (considerably simpler) analogues of the string,
	dilaton, and divisor equations:
\begin{align}
\Brac{\mtau_{a_1}(\gamma_1)\cdots 
\mtau_{a_n}(\gamma_n)\mtau_{0}(T_0)}\!_{g,\beta}^X & =  0 
\label{string}\\[8pt]
\Brac{\mtau_{a_1}(\gamma_1)\cdots 
\mtau_{a_n}(\gamma_n)\mtau_{1}(T_0)}\!_{g,\beta}^X & =  (2g-2) \!\cdot\! 
\Brac{\mtau_{a_1}(\gamma_1)\cdots 
\mtau_{a_n}(\gamma_n)}\!_{g,\beta}^X \label{dilaton}\\[8pt]
\Brac{\mtau_{a_1}(\gamma_1)\cdots 
\mtau_{a_n}(\gamma_n)\mtau_{0}(D)}\!_{g,\beta}^X & = \textstyle{\int_\beta 
D} \!\cdot\! \Brac{\mtau_{a_1}(\gamma_1)\cdots 
\mtau_{a_n}(\gamma_n)}\!_{g,\beta}^X \label{divisor}
\end{align}
where in the last equation (the divisor equation), $D \in
H^2(X,\Q)$.  Exactly as for the tautological psi integrals, there are two
special cases for $(g,\beta)=(1,0)$, namely 
\begin{equation}\label{specialdilaton}
\brac{\mtau_1(T_0)}\!_{1,0}^X = \brok{1}{24} \chi(X) 
\quad \text{ and } \quad 
\brac{\mtau_0(D)}\!_{1,0}^X = - \brok{1}{24}\int_X D\cp c(T_X).
\end{equation}
(Observation: The symbols $\mtau_{a}(\gamma)$ are not defined in the case
$(g,\beta)=(0,0)$, but if $\mtau_0$ were defined to be the usual primary fields
(evaluation classes only), then there would be the usual two special cases
$\brac{\mtau_0(\gamma_1) \mtau_0(\gamma_2)
\mtau_0(1)}\!_{0,0}^X = \int_X \gamma_1 \cp \gamma_2$
and $\brac{\mtau_0(\gamma_1) \mtau_0(\gamma_2)
\mtau_0(D)}\!_{0,0}^X = \int_X \gamma_1 \cp \gamma_2\cp D$.)
%
%
%
%
%
%
\end{blanko}

\subsection{Diagonal classes and the splitting formula}
\label{Sec:diagonal}

\begin{blanko}{Diagonal classes.}\label{diagonal}
	Assume $g>0$ or $\beta>0$.
	The $ij$'th diagonal class $\deltaclass_{ij}\in H^2(\M_{g,S}(X,\beta),\Q)$ is
	defined as the pull-back from $\M_{g,\{p_i,p_j\}}(X,\beta)$ of the Cartier
	divisor $D_{ij}$.  Clearly it is invariant under pull-back along forgetful
	morphisms.  It can also be described as the sum of all boundary divisors
	having $p_i$ and $p_j$ together on a contracted rational tail.  Similarly,
	let $\deltaclass_{ijk}$ denote the sum of all boundary
	divisor having $p_i$, $p_j,$ and $p_k$ together on a contracted rational
	tail.  Note that this class is zero integrated against classes which are
	compatible with pull-back along forgetful morphisms.  This follows from the
	projection formula and the fact that its push-down is zero (no components of
	the curves in this divisor are destabilized).
\end{blanko}
\begin{lemma}\label{push-delta}
	We have $\pi_i{}\lowerstar \deltaclass_{ij} = 1$.  Hence,
	$\pi_i{}\lowerstar \big( \deltaclass_{ij} \cap [\M_{g,S}(X,\beta)]\virtual
	\big) = [\M_{g,S\shortsetminus\{p_i\}}(X,\beta)]\virtual .
	$
\end{lemma}
\begin{dem}
	$\deltaclass_{ij}$ is a sum of boundary divisors.  One of them is $D_{ij} =
	\sigma_j{}\lowerstar (1)$, and being a section, its push-down is the 
	fundamental class.  The other components of $\deltaclass_{ij}$ parameterize
	maps with at least four special points on the contracted tail, which
	therefore are not destabilized when forgetting $p_i$.  Hence these components
	of $\deltaclass_{ij}$ drop dimension and have zero push-down.  (This implies
	in particular that the push-down of $\deltaclass_{ijk}$ is zero.)
The statement about virtual fundamental classes follows from the first 
statement together with the fact that the virtual fundamental class of
the pointed space is the flat pullback of the virtual fundamental class
on the unpointed space.
\end{dem}

\begin{lemma}\label{swap}
	The diagonal classes enjoy the following properties.
	\begin{eqnarray*}
		\deltaclass_{ij}\, \evalmap_i\upperstar (\gamma) 
		& = & \deltaclass_{ij}\, \evalmap_j\upperstar (\gamma),  
		\qquad \text{ for } \gamma\in H\upperstar (X,\Q)\\
		-\deltaclass_{ij}^2 \ 
		= \ \deltaclass_{ij}\,\mpsiclass_i  & = & 
		\deltaclass_{ij}\,\mpsiclass_j\\
		\deltaclass_{ij}\,\deltaclass_{ik} & = 
		& \deltaclass_{ij}\,\deltaclass_{jk} .
	\end{eqnarray*}%
\end{lemma}
These properties motivate the name ``diagonal class''.

\begin{dem}
	Since the involved classes are compatible with pull-back under forgetful
	maps, it is enough to prove the formulas for the two-pointed space
	$\M_{g,\{p_i,p_j\}}(X,\beta)$ (and for the three-pointed space for the third
	formula).  The first formula reads $D_{ij} \, \evalmap_i\upperstar (\gamma) =
	D_{ij} \, \evalmap_j\upperstar (\gamma)$.  This identity follows directly from
	the fact that on the space $\M_{0,\{p_i,p_j,x\}}(X,0)$ corresponding to the
	degree $0$ factor of the divisor $D_{ij}$ ($x$ is the gluing mark), the two
	evaluation morphisms $\nu_i$ and $\nu_j$ coincide.
	
	The second formula follows from the self-intersection formula $D_{ij}^2 = -
	D_{ij} \cdot\pi_j\upperstar \psiclass_i$, and the observation that on the
	two-pointed space $\M_{g,\{p_i,p_j\}}(X,\beta)$ we have $\pi_j\upperstar
	\psiclass_i = \mpsiclass_i$.
	
	The third formula amounts to the following equivalence of intersections of 
	Cartier divisors (on the three-pointed space):
	$$
	(D_{ij} + D_{ijk})(D_{ik}+D_{ijk}) = (D_{ij} + D_{ijk})(D_{jk}+D_{ijk}) .
	$$
	Here the only non-trivial part is $D_{ijk} D_{ik} = D_{ijk} D_{jk}$.  
These
	are transversal intersections, so the formula follows from the fundamental
	linear equivalence $(p_i, p_k\mid p_j,x)= (p_j,p_k\mid p_i,x)$ in the space
	$\M_{0,\{p_i,p_j,p_k,x\}}(X,0)$ corresponding to the degree $0$ factor of the
	divisor $D_{ijk}$ (again $x$ is the gluing mark).
\end{dem}

\begin{BM}\label{D_ij}
	The swapping properties of the lemma hold also with $D_{ij}$ in place of
	$\deltaclass_{ij}$.
\end{BM}

\begin{cor}
	Any integral involving diagonal classes, modified psi classes and
	evaluation classes can be expressed as one involving only modified psi
	classes and evaluation classes.
\end{cor}
\begin{dem}
If there is a diagonal class, say $\deltaclass_{ij}$, use the self-intersection
formula to reduce to the case where there is only one factor of
$\deltaclass_{ij}$.  Next use Lemma~\ref{swap} to substitute all
other occurrences of the index $i$ by $j$.  Now all other classes in the product
are pull-backs from the space without mark $p_i$, so pushing down forgetting
$p_i$ we get rid of $\deltaclass_{ij}$ via the projection formula and
Lemma~\ref{push-delta}.
\end{dem}

For $g\geq 1$, let $N$ be the substack in $\M_{g-1,S\cup\{x',x''\}}(X,\beta)$
consisting of maps $\mu$ such that $\mu(x')=\mu(x'')$.  Its virtual class is given
as the Gysin pull-back $\Delta^{!}[\M_{g-1,S\cup\{x',x''\}}(X,\beta)]\virtual$
in the fiber square
\begin{diagram}[w=9ex,h=4.5ex,tight]
N		& \rInto^{\jmath_N}		& \M_{g-1,S\cup\{x',x''\}}(X,\beta)		\\
\dTo		&  		& \dTo>{(\evalmap_{x'},\evalmap_{x''})}		\\
X		& \rInto_\Delta		& X\times X
\end{diagram}
Let $D_N$ denote the virtual boundary divisor which is the image of the
clutching morphism $\rho_N : N \pil \M_{g,S}(X,\beta)$ which glues the two marks
$x'$ and $x''$ getting a nodal source curve.
%
\begin{lemma}\label{DN}
	With the morphisms $\M_{g,S}(X,\beta) \lInto^{\rho_N} N 
\rInto^{\jmath_N} \M_{g-1,S\cup\{x',x''\}}(X,\beta)$ defined as above, we have
	$$
	\rho_N\upperstar \mpsiclass_i 
	= \jmath_N\upperstar \big(\mpsiclass_i + \deltaclass_{ix'} + 
	\deltaclass_{ix''} - \deltaclass_{ix'x''}\big).
	$$
	
\end{lemma}
\begin{dem}
	Since all the involved classes are compatible with pull-back along forgetful
	morphisms, it is enough to prove the formula for $S=\{p_i\}$.  On
	$\M_{g,\{p_i\}}(X,\beta)$, we have $\mpsiclass_i=\psiclass_i$, and it is
	well-known that $\rho_N\upperstar \psiclass_i = \jmath_N\upperstar
	\psiclass_i$.  It remains to notice that on
	$\M_{g-1,\{p_i,x',x''\}}(X,\beta)$ we have $\psiclass_i = \mpsiclass_i +
	\deltaclass_{ix'} + \deltaclass_{ix''} - \deltaclass_{ix'x''}$, by definition
	of the modified psi class and the pull-back formula for $\psiclass_i$.
\end{dem}

For each stable triple of partitions $S'\cup S'' = S$, $g'+g''=g$, and
$\beta'+\beta''=\beta$ there is a (virtual) boundary divisor $D=D(S',g',\beta'
\mid S'', g'', \beta'')$ defined as the image of the gluing morphism
$$
\M_{g',S'\cup \{x'\}}(X,\beta') 
\times_X \M_{g'',S''\cup \{x''\}}(X,\beta'')  \stackrel{\rho_D}{\langpil}
\M_{g,S}(X,\beta),
$$
with virtual class induced from the virtual classes of the factors.  Precisely,
if $\jmath_D : \M' \times_X \M'' \pil \M'\times \M''$ is the inclusion of that
fibered product in the cartesian product, then $D = \rho_D{}\lowerstar \Delta^! 
\big( [\M']\virtual\boxtimes [\M'']\virtual)$, where $\Delta: X \pil X\times X$
is the diagonal embedding.

\begin{lemma}\label{restr-mpsi}
For a boundary divisor $D=D(S',g',\beta'\mid S'',g'',\beta'')$ with
$(g',\beta')\neq(0,0)$ and $(g'',\beta'')\neq (0,0)$, and with notation as 
above, we have
	$$
	\rho_D\upperstar \mpsiclass_i = \jmath_D\upperstar \big(\mpsiclass_i 
	+\deltaclass_{ix'}\big),
	$$
	assuming $p_i \in S'$. (The mark $x'$ is the gluing mark of that component.)
\end{lemma}

\begin{dem}
	Similar to the proof of \ref{DN}.
\end{dem}
%

To state the following splitting lemma, some further notation is needed. 
Observe that the diagonal classes appearing as correction terms in
\ref{restr-mpsi} are all of the form $\deltaclass_{ix'}$, including the gluing
mark $x'$ ($\deltaclass_{ix''}$ on the other component).  The effect of
multiplying with such a class is to move all classes of mark $p_i$ over to the
gluing mark $x'$, so let us introduce a shorthand notation for this.  Given a
set of marks $B$, and a product of classes $ \smallprod{i\in B}{}
\mtau_{m_i}(\gamma_i)$ indexed by $B$, set
$$
\gamma^B \df \smallprod{p_i\in B}{} \gamma_i, \qquad \text{ and } \qquad
\fat{m}_B \df \sum_{p_i\in B} (m_i -1).
$$
In the definition of the integer $\fat{m}_B$, a priori there might be
negative summands, but we are going to preclude this in the application.
With this notation,

\begin{lemma}\label{split}\textsc{Splitting lemma.}
Let $D= D(S',g',\beta' | S'',g'',\beta'')$ be a boundary divisor with
$(g',\beta')\neq(0,0)$ and $(g'',\beta'')\neq (0,0)$.  Suppose we are given an
almost-top product $\mtau_{m_1}(\gamma_1) \cdots \mtau_{m_n}(\gamma_n)$, then
the integral $\brac{ D \cdot \mtau_{m_1}(\gamma_1) \cdots
\mtau_{m_n}(\gamma_n)}\!_{g,\beta}^X$ is equal to
	$$
	\sum_{e,f} \sum
	\bigg( \Brac{ \smallprod{p_i \in A'}{} \mtau_{m_i}(\gamma_i) \cdot
	\mtau_{\fat{m}_{B'}}(\gamma^{B'}\cp T_e) }\!_{g',\beta'}^X \bigg)
	g^{ef}
		\bigg( \Brac{ \smallprod{p_i \in A''}{} \mtau_{m_i}(\gamma_i) \cdot
	\mtau_{\fat{m}_{B''}}(\gamma^{B''}\cp T_f) }\!_{g'',\beta''}^X \bigg)
	$$
	{\rm The outer sum is over the splitting indices $e$ and $f$ running from $0$
	to $r$.  The inner sum is over all partitions $A'\cup B' = S'$ such that
	$m_i>0$ for all $p_i\in B'$, and over all partitions $A''\cup B'' = S''$ such
	that $m_i>0$ for all $p_i\in B''$.  }
\end{lemma}

\begin{dem}
The splitting indices come from the fact that the virtual class of $D$ is
$\Delta^![\M' \times \M'']\virtual$, the Gysin pull-back of the
diagonal in $X\times X$ along the product of the evaluation morphisms at the two
gluing marks.  This explains the factor
$$
\sum_{e,f} \evalmap_{x'}\upperstar 
(T_e) \, g^{ef} \, \evalmap_{x''}\upperstar (T_f).
$$
Now the cohomology classes at the gluing marks are not left at that, because of
the appearance of the diagonal classes when the modified psi classes restrict to
$D$.  Let us explain what happens on the one-primed component.  Given $p_i\in
S'$, the factor $\mtau_{m_i}(\gamma_i) = \mpsiclass_i^{m_i} \cp
\evalmap_i\upperstar (\gamma_i)$ restricts to $D$ giving 
$$
(\mpsiclass_i +
\deltaclass_{ix'})^{m_i} \cp \evalmap_i\upperstar (\gamma_i) 
= \mpsiclass_i^{m_i} \cp
\evalmap_i\upperstar (\gamma_i) \ 
+ \ \mpsiclass_i^{m_i-1}\cp\deltaclass_{ix'} \cp
\evalmap_i\upperstar (\gamma_i),
$$
after expanding and repeated use of the self-intersection formula in \ref{swap}. 
For each $p_i\in S'$ there is such a sum; expanding the product of all these
sums, we get a sum over all 2-partitions $A'\cup B'=S'$, where the $A'$-part
corresponds to taking the terms without a diagonal class involved, and the
$B'$-part corresponds to the terms with a diagonal class.  Now only at marks
$p_i$ where $m_i>0$ does a diagonal class appear at all, so we must consider
only partitions such that $p_i\in A'$ whenever $m_i=0$.  Now the effect of the
diagonal class $\deltaclass_{ix'}$ is to move the classes
$\mpsiclass_i^{m_i-1}\evalmap_i\upperstar (\gamma_i)$ over to the gluing mark
$x'$.  To be more precise, use \ref{swap} to replace the indices $i$ by $x'$,
and then push down forgetting $p_i$ (cf.~\ref{push-delta}).  Doing this for each
of the marks in $B'$, accumulates at $x'$ the class
$\mtau_{\fat{m}_B}(\gamma^B)$, and there was already a class
$\evalmap_{x'}\upperstar (T_e)$, thus totalling $\mtau_{\fat{m}_B}(\gamma^B\cp
T_e)$, as claimed.  This explains the contribution from the one-primed component. 
The same arguments hold for the other component.
\end{dem}

\section{Topological recursion relations in genus $0$ and $1$}
\label{Sec:toprec}
\setcounter{lemma}{0}

The topological recursion relations for modified psi classes rely on the same
two facts that drive the topological recursion for the usual descendants: that
the psi classes admit an expression in terms of boundary divisors, and that it
is known how to restrict to such boundary divisors (the splitting lemma). 
However there is a crucial difference, namely that the modified psi classes do
not restrict to boundary divisors in the straightforward way the tautological
classes do, so the corresponding topological recursions are much more
complicated.  For first descendants however, which are the ones needed to
describe tangency conditions, the deformed metric of~\cite{Kock:0006} allows
for a concise way of writing the equations.

\begin{blanko}
{Modified psi classes as boundary corrections.} In any case, the modified psi
class is a boundary correction to the tautological psi class.  Precisely, on 
any stack $\M_{g,S}(X,\beta)$ with
$(g,\beta)\neq (0,0)$, let $\boundary_i$ denotes the sum of all boundary
divisors having $p_i$ on a contracted rational tail.  Then
	$$
	\mpsiclass_i = \psiclass_i - \boundary_i.
	$$
This follows readily from the fact that $\boundary_i$ is zero for $n=1$ (by 
stability), and that it pulls back along forgetful morphisms in the same manner 
as $\psiclass_i$, to wit: $\pi_0\upperstar \boundary_i = \boundary_i - D_{i,0}$.
\end{blanko}

\subsection{Genus zero}
  
Recall that in genus zero, and when there are at least three marks,  the
tautological psi class admits an expression in terms of boundary divisors;
precisely,
$$
\psiclass_i = (p_i \mid p_j,p_k),
$$
for any two fixed marks $p_j$ and $p_k$, distinct, and distinct from $p_i$. 
(See e.g.\ Getzler~\cite{Getzler:9801}.)  Therefore also the modified psi classes
can be written in terms of boundary divisors.

\begin{satz}\textsc{Topological recursion (genus zero).}\label{rec}
	The following recursive relation holds on a space with at least three marks,
	say $\M_{0, S\cup\{p_1,p_2,p_3\}}(X,\beta)$.  
	\small
	\begin{multline*}
	\Brac{ \mtau_{m_1+1}(\gamma_1) 
	\mtau_{m_2}(\gamma_2)\mtau_{m_3}(\gamma_3)\cdot
	\smallprod{p_i \in S}{}\mtau_{m_i}(\gamma_i)}\!_\beta
	=\\ + 
	\Brac{ \mtau_{m_1}(\gamma_1) \mtau_{m_2+m_3}(\gamma_2\cp \gamma_3)
	\cdot \smallprod{p_i \in S}{}\mtau_{m_i}(\gamma_i)}\!_\beta \\
	-
	\Brac{ \mtau_{m_1+m_2}(\gamma_1\cp\gamma_2) \mtau_{m_3}(\gamma_3)
	\cdot \smallprod{p_i \in S}{}\mtau_{m_i}(\gamma_i)}\!_\beta \\
	-
	\Brac{ \mtau_{m_1+m_3}(\gamma_1\cp\gamma_3) \mtau_{m_2}(\gamma_2)
	\cdot \smallprod{p_i \in S}{}\mtau_{m_i}(\gamma_i)}\!_\beta \\
	+ \bigsum{}{} \sum
	\big(  \Brac{\smallprod{i\in 
	A'}{}\mtau_{m_i}(\gamma_i) \cdot \mtau_{\fat{m}_{B'}}
	(\gamma^{B'}\cp T^e)}\!_{\beta'}\big) 
	g^{ef} 
\big(  \Brac{\smallprod{i\in 
	A''}{}\mtau_{m_i}(\gamma_i) \cdot \mtau_{\fat{m}_{B''}}
	(\gamma^{B''}\cp T^f)}\!_{\beta''}\big) 	
\end{multline*}
\normalsize The big outer sum is over all stable splittings $S'\cup S'' = S$,
$\beta'+\beta''=\beta$, $\beta'>0$, $\beta''>0$, and over the splitting indices
$e,f$ running from $0$ to $r$.  The inner sum is over all partitions $A'\cup B'
= S'\cup \{p_1\}$ such that $m_i>0$ for all $i\in B'$, and over all partitions
$A''\cup B''
= S''\cup \{p_2,p_3\}$ such that $m_i>0$ for all $i\in B''$.
\end{satz}

\begin{dem}
	Write the first modified psi class as a sum of boundary divisors,
	\begin{equation}
	\mpsiclass_1 = (p_1 \mid p_2,p_3) - \boundary_1,\label{mpsi=boundary}
	\end{equation}
	and compute the restriction of the remaining factors of the product to each
	of the irreducible components of this expression.  Let us first consider the
	boundary divisors involving a contracted tail (called \soft boundary
	divisors).  Since there are at least two marks on the contracted tail, say
	$p_i$ and $p_j$, we are in position to use the Remark~\ref{D_ij}: all the
	remaining factors are invariant under pull-back so we can substitute all
	indices $i$ by $j$, and then compute the integral via a push-down, forgetting
	$p_i$.  Now, unless $p_i$ and $p_j$ are the only marks on the 0-tail, we get
	zero push-down, cf.~\ref{push-delta}.  So among the \soft boundary divisors in
	(\ref{mpsi=boundary}), only the ones with just two marks on the
	contracted tail contribute.  From $(p_1\mid p_2,p_3)$ we get $D_{23}$ as
	well as $D_{1k}$ for $p_k\in S$.  Now all the latter are also in
	$\boundary_1$.  On the other hand, in $\boundary_1$ we find $D_{12}$ and
	$D_{13}$ which are not in $(p_1|p_2,p_3)$.  So all together only three \soft
	components contribute:
$$
D_{23} - D_{12}-D_{13}.
$$
The effect of restricting to such a boundary divisor, say $D_{12}$ is to merge
all the classes indexed by $1$ and $2$, so this explains the first three terms 
in the formula.

Now for the hard boundary: We are concerned with the sum of all boundary
divisors $D= D(S'\cup \{p_1\}, \beta' \mid S''\cup \{p_2, p_3\}, \beta'')$
with $\beta'>0$ and $\beta''>0$.  (Thus explaining this summation in the
formula.)  For each of these boundary divisors the splitting lemma applies.
\end{dem}

\begin{BM}
	This recursion relation determines all enumerative descendants from the
	Gromov-Witten invariants, and while it is not very pleasant to look at, it is
	quite effective: each step reduces the number of modified psi
	classes by one, so the depth of recursion is equal to the number of modified
	psi classes.
\end{BM}

\subsection{First enumerative descendants}

For the sake of describing tangency conditions, only invariants with exponent at
most one on each psi class are needed (cf.~Section~\ref{Sec:char}).  
These invariants are very well-behaved, and the recursion can be written in
a nice manner, as we now proceed to describe, following \cite{Kock:0006}.

\begin{blanko}{The tangency quantum potential.}\label{Gamma}
Set
$$
\brac{\mbtau_0^{\fat{a}} \; \mbtau_1^{\fat{b}}}\!_{g,\beta}^X \df 
\Brac{ \smallprod{k=0}{r} 
(\mtau_0(T_k))^{a_k}(\mtau_1(T_k))^{b_k}}\!_{g,\beta}^X,
$$
where $\fat{a}=(a_0,\ldots,a_r)$ and $\fat{b}=(b_0,\ldots,b_r)$ are vectors of
non-negative integers.  The integral is zero unless $\sum_i \deg(T_i)
(a_i+b_i+1) = \vdim \M_{g,n}(X,\beta)$, where $n= \sum (a_i + b_i)$ is the 
number of marks.
The generating function corresponding to these first enumerative descendants
is called the {\em tangency quantum potential}:
\begin{align*}
\Gamma^{(g)}(\fat{x},\fat{y}) &=
\sum_{\beta} 
\brac{\exp \big( \smallsum{}{} x_i \mtau_0(T_i) 
+ \smallsum{}{} y_i \mtau_1(T_i)\big)}\!_{g,\beta}^X \\
&= \sum_{\beta}\sum_{\fat{a},\fat{b}} 
\frac{\fat{x}^{\fat{a}}}{\fat{a}!}\frac{\fat{y}^{\fat{b}}}{\fat{b}!}
\brac{\mbtau_0^{\fat{a}} \; \mbtau_1^{\fat{b}}}\!_{g,\beta}^X
\end{align*}
where $\fat{x} = (x_0,x_1,\ldots,x_r)$ and $\fat{y}=(y_0,y_1,\ldots,y_r)$ are
formal parameters, and   we employ the usual multi-index notation $\fat{s}!  =
s_0!  s_1!  \cdots s_r!$, and  $\fat{y}^{\fat{s}}=
y_0^{s_0}y_1^{s_1}\cdots y_r^{s_r}$.

For $g=0$, the degree expansion is only over $\beta>0$.
%
The variables $\fat{x}$ are the usual formal variables from quantum cohomology,
so when $\fat{y}$ is set to zero, $\Gamma^{(0)}$ reduces to the usual (quantum
part of the) genus zero Gromov-Witten potential.

Note that, for $g=1$, the degree $0$ part of the tangency quantum potential
consists only of a couple of terms, corresponding to the special cases of the
dilaton and divisor equations~(\ref{specialdilaton}).  In higher genus, the only
other occurrences of a non-zero degree $0$ part are for $g=2$, when $X$ is a
surface.
\end{blanko}

\begin{blanko}{Derivatives.}\label{der}
	First, observe that $\Gamma_{x_k} \df \frac{\d}{\d x_k}\Gamma$ is the
	generating function for the numbers
	$\brac{\mtau_0(T_k)\;\mbtau_0^{\fat{a}}\,\mbtau_1^{\fat{b}} }\!_{g,\beta}^X$,
	in the sense that this number is the coefficient of
	$\frac{\fat{x}^{\fat{a}}}{\fat{a}!}\frac{\fat{y}^{\fat{b}}}{\fat{b}!}$ in
	$\Gamma_{x_k}$.  Next, adopt the convention that if two $x$-variables are
	coupled in a parenthesis then the meaning is
\begin{equation}\label{Gammaij}
	\Gamma_{(x_i x_j)} \df \sum_{k=0}^r g_{ij}^k \ \Gamma_{x_k}.
\end{equation}
	This function is the ``directional derivative with respect to the product
	\newcommand{\rightdef}{={\raisebox{0.255ex}{\normalfont\scriptsize\!\!:}}\: }
	$T_i\cp T_j \rightdef \sum g_{ij}^k T_k$''.  Precisely, $\Gamma_{(x_i x_j)}$
	is the generating function for the invariants $\brac{\mtau_0(T_i\cp
	T_j)\;\mbtau_0^{\fat{a}}\,\mbtau_1^{\fat{b}} }\!_{g,\beta}^X$ .
\end{blanko}

\begin{blanko}
	{The deformed metric.} (See Kock~\cite{Kock:0006} for details.)  Consider the
	non-degenerate symmetric bilinear pairing $H\upperstar (X,\Q) \pil \Q[[\fat
	y]]$ given by the tensor elements
	$$
	\gamma_{ij} \df \sum_{\fat{s}}
	\frac{(-2\fat{y})^{\fat{s}}}{\fat{s}!} \int_X \fat{T}^{\fat{s}} \cp T_i \cp
	T_j .
	$$
	The inverse matrix of $(\gamma_{ij})$ is given by
	$$
	\gamma^{ij}=
	\sum_{\fat{s}}
	\frac{(2\fat{y})^{\fat{s}}}{\fat{s}!} 
	\int_X \fat{T}^{\fat{s}} \cp\check T_i \cp
	\check T_j ,
	$$
	where $\check T_k = \sum_{m} g^{km} T_m$.
	When the formal parameters of $\gamma$ are reset, we recover the
	corresponding $g$, for example $\gamma^{ij}(\fat{0}) = g^{ij}$.
\end{blanko}

This metric encodes the combinatorics of the diagonal classes
that appear in the splitting lemma, and in this formalism, the topological
recursion relation for the genus $0$ first enumerative descendants takes the
following simple form.

\begin{satz}\label{PDE}
	The genus $0$ tangency potential satisfies the differential equations
	$$
	\Gamma^{(0)}_{y_k x_i x_j } = \Gamma^{(0)}_{x_k(x_i x_j)} 
	- \Gamma^{(0)}_{(x_k x_i) x_j}
	- \Gamma^{(0)}_{(x_k x_j) x_i} +
	\sum_{e,f} \; \Gamma^{(0)}_{x_k x_e} \ 
	\gamma^{ef} \ \Gamma^{(0)}_{x_f x_i x_j}.
	$$
\end{satz}

\begin{blanko}{Integration of the differential equation.}\label{integration}
In the particular differential equation 
where $k=i=j$,
all the terms are total derivatives with respect to $x_k$, so
we actually get a better differential equation
$$
\Gamma^{(0)}_{x_k y_k} = - \Gamma^{(0)}_{(x_k x_k)}
+
\brok{1}{2} \sum_{e,f} \; \Gamma^{(0)}_{x_k x_e} 
\ \gamma^{ef} \ \Gamma^{(0)}_{x_f x_k}.
$$
(Since all terms are exponential in $x_k$, the integration constant must be 
zero.)	
\end{blanko}
\begin{eks}\label{gammaP2}
	For $\P^2$ (with $h\df c_1(\OO(1))$ and basis $T_i = h^i$), we have
	$$
	(\gamma^{ij}) = \exp(2y_0)\begin{pmatrix}
	\phantom{0}0\phantom{0} & 
	\phantom{0}0\phantom{0}& \phantom{0}1\phantom{0} \\
	0 & 1 & 2y_1 \\
	1 & 2y_1 & 2y_1^2 + 2y_2
	\end{pmatrix} .
	$$
With $k=1$, the differential equation reads
$$
\Gamma^{(0)}_{y_1 x_1} =  - \Gamma^{(0)}_{x_2}
+	\brok{1}{2}\Gamma^{(0)}_{x_1 x_1}\,\Gamma^{(0)}_{x_1 x_1}
+	2y_1 \Gamma^{(0)}_{x_1 x_1} \, \Gamma^{(0)}_{x_1 x_2}
+	(y_1^2 + y_2)\Gamma^{(0)}_{x_1 x_2}\,\Gamma^{(0)}_{x_1 x_2}.
$$
\end{eks}

\subsection{Genus one}

In genus one, there is the following expression for the tautological psi class,
(see~\cite{Getzler:9801}).
$$
\psiclass_i = \brok{1}{24} D_N + B_i ,
$$
where $B_i$ is the sum of all boundary divisors such that $p_i$ is on a rational
component.  It follows that $\mpsiclass_i$ is expressed as $\brok{1}{24}D_N$
plus the sum of all boundary divisor such that $p_i$ is on a {\em
non-contracted} rational component.

Hence there is a topological recursion relation for the genus $1$ enumerative
descendants, analogous to \ref{rec}, as it easily follows from the same
arguments.  We state only the formula for the first enumerative descendants:

\begin{prop}\label{TRR1}\textsc{Topological recursion (genus one).}
	The genus $1$ tangency quantum potential $\Gamma^{(1)}$ satisfies the
	differential equations
	$$
\Gamma^{(1)}_{y_k} = \sum_{e,f} \Gamma^{(0)}_{x_k x_e} \; \gamma^{ef} \; 
\Gamma^{(1)}_{x_f} + 
\frac{1}{24} \sum_{e,f} \gamma^{ef} \Gamma^{(0)}_{x_k x_e x_f} .
$$
\end{prop}
For $k=0$, the right hand side needs an extra term $+\frac{1}{24}\chi(X)$ to
account for the special case of the dilaton equation.  This is necessary since
$\Gamma^{(0)}$ is not defined in degree $0$.
Otherwise, 
note that this is exactly like the corresponding topological recursion relation 
for the tautological psi classes, only with the deformed metric in place of the 
Poincar\'e metric.

\section{Tangency conditions and characteristic numbers}
\label{Sec:char}
\setcounter{lemma}{0}

The motivation for studying modified psi classes comes from enumerative
geometry: we show that tangency conditions are naturally expressed
in terms of modified psi classes, and that characteristic numbers (of rational
curves) are first enumerative descendants.  

Classically, characteristic numbers (of $\P^2$) are the numbers of curves of
given genus and degree which pass through given points and are tangent to given
lines.  The generalization to other homogeneous varieties involves the choice of
what subvarieties one should impose conditions with respect to.  In the context
of first enumerative descendants, the natural conditions to include are those of
being tangent to a hypersurface, or more generally: being tangent to a 
hypersurface along a specified subvariety.  (It should be
stressed that the inclusion of these compound conditions is not necessary for
the recursions to work.)

The approach is a natural extension of the way  Gromov-Witten invariants
describe incidence conditions.
%

\bigskip

From now on, $X$ will denote a homogeneous variety.  In this section, for a
homology cycle or subvariety $W$, abusively we let $[W]$ denote the cohomology
class representing $W$.

\begin{blanko}{Incidence conditions and Gromov-Witten invariants.}
	Put $\M_{\emptyset} = \M_{0,\emptyset}(X,\beta)$, let $W\subset X$ be a
	subvariety of codimension $c$, and let $\incidencelocus{W}\subset
	\M_{\emptyset}$ denote the locus of maps whose image is incident to $W$. 
	Under suitable circumstances, $\incidencelocus{W}$ is of codimension $c-1$. 
	The goal is to compute the intersection $\incidencelocus{W_1} \cap \cdots
	\cap \incidencelocus{W_n}$ for general $W_1,\ldots,W_n$ whose codimensions
	add up to $\dim \M_{\emptyset} + n$.  The solution is a Gromov-Witten
	invariant: Let $\pointedincidencelocus{p_i}{W_i} \subset
	\M_{\{p_i\}}=\M_{0,\{p_i\}}(X,\beta)$ denote the inverse image stack of $W_i$
	under the evaluation morphism $\evalmap_i : \M_{\{p_i\}} \pil X$, i.e., the
	locus of maps whose mark lands in $W_i$.  For general $W_i$ of codimension
	$c_i$, $\pointedincidencelocus{p_i}{W_i}$ is reduced of codimension $c_i$ and
	is represented by the evaluation class $\evalmap_i\upperstar [W_i]\in
	H\upperstar (\M_{\{p_i\}}, \Q)$.  If $c_i\geq 2$, and the general map in
	$\M_{\emptyset}$ is generically 1--1, then for general $W_i$, the restriction
	of the forgetful morphism $\pointedincidencelocus{p_i}{W_i} \pil
	\incidencelocus{W_i}$ is also generically 1--1.  In this case,
	$[\incidencelocus{W_i}] = \pi_i{}\lowerstar \evalmap_i\upperstar [W_i]$.
		
%
%
		
	Now a crucial property of the evaluation classes is their 
	compatibility with forgetful morphisms. Concretely, if $\hat\pi_i :
	\M_S \pil \M_{\{p_i\}}$ is the forgetful morphism that forgets all marks in 
	$S=\{p_1,\ldots,p_n\}$ except $p_i$ then
	$$
	\hat\pi_i\upperstar \evalmap_i\upperstar [W_i] =
	\evalmap_i\upperstar [W_i] \quad \text{ in } H\upperstar (\M_S,\Q).
	$$
	So we can pull the various pointed classes back to $\M_S$ and integrate here. 
	By the projection formula, we see that
	$$
	\int_{\M_{\emptyset}} [\incidencelocus{W_1}] \cp \cdots \cp
	[\incidencelocus{W_n}] = \int_{\M_S} \evalmap_1\upperstar [W_1] \cp
	\cdots \cp \evalmap_n\upperstar [W_n]
	$$
	--- this last number is a Gromov-Witten invariant.
%
%
\end{blanko}
%
%
%
%
%
%
%
%
%
%
%

The same viewpoint will now be applied to tangency conditions.  The tangency
condition is described in terms of modified psi classes by going up to a
one-pointed space.  By construction, modified psi classes are compatible with
forgetful morphisms, so the integral of these markless classes can also be
computed on the $n$-pointed space where it is a first enumerative descendant.


\subsection{Tangency conditions via modified psi classes}

Let now $\M_{\emptyset}$ denote a substack in $\M_{g,\emptyset}(X,\beta)$ (e.g.,
the whole stack), and let $\M_{\{p_i\}}$ be its stack-theoretic inverse image
under the forgetful morphism $\pi_i : \M_{g,\{p_i\}}(X,\beta) \pil
\M_{g,\emptyset}(X,\beta)$.

\begin{blanko}{Virtual tangency conditions.}
	Let $V\subset X$ be a hypersurface, and let $N_V$ denote its normal bundle.
The virtual tangency condition to $V$ (at the marked point $p_i$) is defined
by the following method.  The differential gives a natural map of bundles on
$\M_{\{p_i\}}$:
\begin{equation}
\label{diffff}
d\evalmap_i: \LL\dual_i \pil \evalmap_i\upperstar (T_X).
\end{equation}
Consider the substack 
$\pointedincidencelocus{p_i}{V}=\evalmap_i^{-1}(V) \subset \M_{\{p_i\}}$
of maps $[\{p_i\}\pil C\stackrel{\mu}{\pil} X]$ for which $\mu(p_i) \in V$.
Since $\evalmap_i$ restricted to $\pointedincidencelocus{p_i}{V}$ factors through
$V$, there is a natural sequence on $\pointedincidencelocus{p_i}{V}$
obtained from (\ref{diffff}) and the normal sequence on $V$:
\begin{equation}\label{comppp}
	\LL\dual_i \stackrel{d\evalmap_i}{\langpil} \evalmap_i\upperstar (T_X|_V) 
	\langpil \evalmap_i\upperstar (N_V) .
\end{equation}
The degeneracy locus of the composition $\LL\dual_i \pil \evalmap_i\upperstar
(N_V)$ is the locus of maps $\mu$ for which the the differential $d\mu_{p_i}$
has image in the tangent space of $V$, i.e., the maps tangent to $V$ exactly at
the mark $p_i$.  The cohomology class representing this degeneracy locus in
$H^2(\pointedincidencelocus{p_i}{V}, \Q)$ is simply $\evalmap_i\upperstar [V]
+ \psiclass_i$.  Hence the class $ \evalmap_i\upperstar [V] \cp (
\evalmap_i\upperstar [V] + \psiclass_i) \in H^4(\M_{\{p_i\}}, \Q) $ is defined
to represent the virtual tangency condition to $V$.

Consider now the corresponding $n$-pointed stack $\M_S$ with its forgetful
morphisms $\hat \pi_i : \M_S \pil \M_{\{p_i\}}$.  The
virtual tangency class in $H\upperstar (\M_S, \Q)$ is defined to be the
$\hat\pi_i$~pull-back of the corresponding class on
$\M_{\{p_i\}}$.  In other words, the virtual condition of
tangency to $V$ (at mark $p_i$) is of class
$$
\evalmap_i\upperstar [V] \cp ( \evalmap_i\upperstar [V] + \mpsiclass_i) 
\in H^4 (\M_S, \Q).
$$
\end{blanko}

In general, this virtual condition is fulfilled by some maps that we would not
normally think of as tangent to $V$.  
Notably, if the marked point is on a component that contracts to a point of $V$
then the map fulfills the virtual tangency condition.

\begin{blanko}{Tangency along a subvariety.}\label{ThetaVW}
	More generally, let $W\subset V$ be a subvariety, and perform the same
	argument in the substack $\pointedincidencelocus{p_i}{W}$ of maps whose mark
	lands in $W$.  Here, the degeneration of the map $\LL\dual_i \pil
	\evalmap_i\upperstar (N_V)$ is the locus of maps for which furthermore the
	differential $d\mu_{p_i}$ has image in the tangent space of $V$.  The
	cohomology class representing this degeneracy locus in $H\upperstar
	(\M_{\{p_i\}}, \Q)$ is
$$
\evalmap_i\upperstar [W] \cp( \evalmap_i\upperstar [V] + \mpsiclass_i),
$$
which is defined to be the virtual conditions of tangency to $V$ along $W$.
\end{blanko}

\begin{blanko}{Pointed and unpointed conditions.}
	Let $\Theta_i \df \evalmap_i\upperstar [W_i] \cp( \evalmap_i\upperstar [V_i] +
	\mpsiclass_i)$ for some $W_i\subset V_i$.  The virtual class of un-pointed
	tangency is defined to be $\tangencyclass_i \df \pi_i{}\lowerstar
	\Theta_i \in H\upperstar
	(\M_{\emptyset}, \Q)$.  Since the classes $\Theta_i$ are compatible with
	forgetful morphisms, the projection formula yields
	$$
	\int_{[\M_{\emptyset}]\virtual} \tangencyclass_1 \cp \cdots \cp
	\tangencyclass_n = \int_{[\M_S]\virtual} \Theta_1 \cp
	\cdots \cp \Theta_n
	$$
	Ideally, (and assuming the codimension of the classes $\Theta_i$ add up to
	$\vdim \M_S$), the first integral is the number of curves tangent to the $n$
	hypersurfaces (along the specified subvarieties).  The second integral is a
	first enumerative descendant.
%
\end{blanko}

\subsection{Virtual characteristic numbers}

\begin{blanko}{Virtual characteristic numbers.}\label{virtualN}
	For simplicity, we consider as incidence conditions only those with respect
	to the basis elements $T_i$.  Assume $T_1,\ldots,T_p$ are the divisor
	classes.  In general, requiring a curve to be incident to a hypersurface is
	an empty condition, so we consider only $i=p+1,\ldots,r$ (codimension $2$ and
	up).  We denote the corresponding conditions $\Omega_i$, and use the same
	symbol for the virtual class in the moduli space:
	$$
	\Omega_i = \mtau_0(T_i) \qquad \text{ for } i=p+1,\ldots,r.
	$$

	Suppose we have chosen $t$ types of conditions, $\Theta_j$, of being tangent
	to certain hypersurfaces (at a specified mark), possibly along specified
	subvarieties.
	The virtual classes of these conditions
	are expressed (cf.~\ref{ThetaVW}) as linear combinations
\begin{eqnarray}\label{rhosigma}
	\Theta_j & = & \sum_{i=0}^r \rho_{ij} \mtau_0(T_i) + 
	\sum_{i=0}^r \sigma_{ij} \mtau_1(T_i), \qquad \text{ for } j=1,\ldots,t.
\end{eqnarray}%

	Now define virtual characteristic numbers 
	\begin{equation}\label{vnumbers}
	\til N^{(g)}_\beta\!(a_{p+1},\ldots,a_r,
	b_1,\ldots,b_t) \df \Brac{ \; \Omega_{p+1}^{a_{p+1}} \cdots \Omega_r^{a_r}
	\; \cdot \; \Theta_1^{b_1} \cdots \Theta_t^{b_t} \; }\!_{g,\beta}^X ,
	\end{equation}
	This is the virtual number of genus $g$ curves of class $\beta$ satisfying
	$a_i$ conditions of type $\Omega_i$, for $i=p+1,\ldots,r$, and further $b_j$
	conditions of type $\Theta_j$, for $j=1,\ldots,t$, (provided the codimensions
	of all the imposed conditions add up to $\vdim \M_{g,n}(X,\beta)$), where 
	$n=\sum a_i + \sum b_j$).
\end{blanko}

\begin{blanko}{Generating functions.}\label{gen-function}
	Introduce partial degrees $d_i \df \int_\beta T_i$, for $i=1,\ldots,p$, and
	define the generating function for the virtual characteristic numbers,
\begin{align}
	\til G^{(g)}(\fat{u},\fat{v}) &\df   	
	\sum_{\beta>0} \exp( \smallsum{i=1}{p}d_i u_i) 
	\underset{b_1,\ldots,b_t}{\sum_{a_{p+1},\ldots,a_r}} 
\frac{u_{p+1}^{a_{p+1}}}{a_{p+1}!}\dots\frac{u_r^{a_r}}{a_r!}
	\cdot 	\frac{v_1^{b_1}}{b_1!}\cdots\frac{v_t^{b_t}}{b_t!}
		\til N^{(g)}_\beta\!(a_{p+1},\ldots,a_r,b_1,\ldots,b_t),
 \notag\\
	 &=  	\sum_{\beta>0} \sum_{\fat{a},\fat{b}} 
\frac{\fat{u}^{\fat{a}}}{\fat{a}!}\frac{\fat{v}^{\fat{b}}}{\fat{b}!}
		d_1^{a_1}\cdots d_p^{a_p} \,
		\til N^{(g)}_\beta\!(a_{p+1},\ldots,a_r,b_1,\ldots,b_t) ,
\end{align}
with multi-index notation as in \ref{Gamma}.
By definition of the virtual characteristic numbers, 
this can also be written
$$
\til G^{(g)}(\fat u, \fat v) = \sum_{\beta>0} \Brac{ \exp \big( 
\smallsum{i,j}{} ( u_i + v_j \rho_{ij}) \mtau_0(T_i) 
+ \smallsum{i,j}{} v_j \sigma_{ij} \mtau_1(T_i) 
\big)}\!_{g,\beta}^X  .
$$
Hence we get:
\end{blanko}

\begin{lemma}\label{GGamma}
	Let $\Gamma^{(g)}_+$ denote the $\beta>0$ part of 
	the tangency quantum potential.
	The series $\til G^{(g)}$ is related to the 
tangency quantum potential by 
	$$
	\til G^{(g)}(\fat{u},\fat{v}) = \Gamma^{(g)}_+(\fat{x},\fat{y}),
	$$
	subject to the change of variables
\begin{alignat*}{4}
&\phantom{qwerty}\quad&x_0 &= 0, \quad  &   
x_i &= u_i  +\smallsum{j=1}{t} \rho_{ij} v_j                 
\quad &    &\text{ for } i=1,\ldots,r \\
&\phantom{qwerty}\quad&y_0 &= 0, \quad  &  y_i 
&= \phantom{u_i + {}}\smallsum{j=1}{t} 
\sigma_{ij} v_j \quad &    &\text{ for } i=1,\ldots,r   ,
\end{alignat*}
which is just the dual to the change of variables relating the conditions
$\Omega_j$ and $\Theta_j$ to the basis $\mbtau_0$ and $\mbtau_1$.
\qed
\end{lemma}

\subsection{Enumerative significance}

In order for the virtual tangency class to express the enumerative geometry
accurately, some additional assumptions are needed.

\begin{blanko}{Assumption on $W\subset V$.}\label{assumptionWV}
	Let $W$ be of codimension $c$ in $X$.  A sufficient generality condition on
	$W\subset V$ is this: $W$ can be translated to intersect any curve
	transversely at a point.  Of $V$ we require: The linear system consisting of
	the members of $\norm{V}$ that contain $W$ separates normal vectors of $W$. 
	(I.e., for each point $x\in W$ and each normal vector $w\in N_x W$, there is
	a $V\supset W$ such that $w \not\in T_x V$.)
\end{blanko}

\begin{blanko}{Assumptions on the family.}\label{assumptionM}
	A markless stable map $\mu:C\pil X$ is said to have a {\em deep cusp} if
	there is a point where the differential vanishes to (at least) second order. 
	For example any map with a contracted component has deep cusp in this sense.

We will impose the following two conditions on the family $\M_\emptyset$:
\begin{list}{}{%
	\setlength{\topsep}{9pt}\setlength{\leftmargin}{54pt}
	\setlength{\rightmargin}{28pt}\setlength{\labelwidth}{20pt} 
	\setlength{\itemsep}{3pt}\setlength{\parsep}{0pt} } 
	
	\item[$\bullet$] $\M_\emptyset$ is reduced of dimension at most $k$ (for some
	$k$).
	
	\item[$\bullet$] The locus of maps with deep cusp is of dimension at
	most $k-1$.
\end{list}
\end{blanko}

Let $\pointedtangencylocus$ denote the degeneracy locus of (\ref{comppp}), and
let $\tangencylocus$ denote its image in $\M_\emptyset$.
\begin{prop}\label{enum}
	With notation and assumptions as above, for general $W\subset V$ we have: (i)
	The degeneracy locus $\pointedtangencylocus\subset \M_{\{p_i\}}$, and thus
	$\tangencylocus = \pi_i(\pointedtangencylocus) \subset \M_\emptyset$, is of
	dimension at most $k-c$.  (ii) The locus of deep cusp in $\tangencylocus$ is
	of dimension at most $k-c-1$.

(iii) Except possibly for those components of $\pointedtangencylocus$
whose general map has a contracted component, $\pointedtangencylocus$ is
reduced, and the morphism $\pointedtangencylocus \pil
\tangencylocus$ is generically 1--1 and \'etale.  In particular, the
maximal dimension part of $\tangencylocus$ is reduced too.
\end{prop}

\begin{dem}
	The degeneracy locus of (\ref{comppp}) can be considered the zero scheme of a
	section of $\evalmap_i\upperstar \OO_X(V) \tensor \LL_i$.  Outside the locus of
	maps whose differential vanishes at $p_i$, this line bundle is generated by
	the space of global sections corresponding to $V$ containing $W$.  Indeed, by
	assumption on $V\supset W$, for each such map $\mu$ there exists a
	hypersurface $V\supset W$ whose tangent space does not contain the tangent
	vector of $\mu(C)$ at $\mu(p_i) \in W$.  So by
	Kleiman-Bertini~(\cite{Kleiman:transversality}, Remark~7), along this locus
	the section vanishes with multiplicity $1$, for general $V\supset W$.
	
	The locus of maps in $\M_{\{p_i\}}$ whose differential vanishes {\em simply}
	at $p_i$ has dimension at most $k$, and requiring further that the mark maps
	to $W$ brings the dimension down to at most $k-c$ (by the generality of $W$),
	and the section vanishes automatically along this locus, with multiplicity
	$1$ by assumption.  The dimension of the corresponding loci of isolated deep
	cusp is one lower, hence at most $k-c-1$, as asserted.
	
	Now for maps with a contracted component.  In $\M_\emptyset$ they occur in
	dimension at most $k-1$, so in $\M_{\{p_i\}}$ they have dimension at most
	$k$.  Contracted components can also arise as a result of stabilizing when
	the mark ``falls on a node''; this type of contracted component also occurs
	in dimension $k$.  Requiring the mark to map to $W$ cuts the dimension down
	to $k-c$, and along this locus the section automatically vanishes.
	
	Concerning the image of these loci in $\tangencylocus$: if the general map of
	a locus retains the contracted component after forgetting the mark, then
	the dimension drops to $k-c-1$, since one dimension of the locus was the
	freedom of the mark moving on the contracted component.
	
	For general $V$, the morphism $\pointedtangencylocus \pil
	\tangencylocus$ is generically 1--1 and \'etale because for a fixed map
	$\mu$, the set of tangent hypersurfaces $V$ that are not {\em simply tangent}
	is of codimension at least one among all the $V$.  This is clear when
	$X=\P^2$ (a plane curve has only a finite number of bitangents and flexes);
	the general case can be obtained from this fact by embedding $X$ in a
	projective space by the linear system $\norm{V}$ and then projecting down to
	a suitable $\P^2$.
\end{dem}

\begin{blanko}{Enumerative significance.}
	Since the substack $\tangencylocus$ thus inherits the properties of the
	original family, the construction can be applied inductively.  By a dimension
	reduction argument, if the general map of $\M_\emptyset$ is irreducible,
	(resp.\ an immersion, resp.\ 1--1), then the general map in $\tangencylocus$
	is again irreducible, (resp.\ an immersion, resp.\ 1--1).  In particular if
	there is only a finite numbers of solutions, they are all irreducible,
	(resp.\ immersions, resp.\ 1--1).
	
	Hence, in this case, a general top intersection of incidence loci and
	tangency loci consists of a finite numbers of reduced points which correspond
	to irreducible maps, and the cardinality of this set is the integral of the
	corresponding pointed conditions over the fundamental class of the marked
	space.

	Clearly, $\M_{0,\emptyset}(X,\beta)$ satisfies the conditions.  On the other
	hand, for $g\geq 1$, the space $\M_{g,\emptyset}(X,\beta)$ does not satisfy
	the conditions, due to contracted tails.
%
\end{blanko}


\subsection{Characteristic numbers of rational curves}

With notation as in \ref{virtualN}, let
$N^{(0)}_\beta\!(a_{p+1},\ldots,a_r,b_1,\ldots,b_t)$ denote the number of
irreducible, rational curves of class $\beta$ satisfying $a_i$ conditions of
type $\Omega_i$, for $i=p+1,\ldots,r$, and further $b_j$ conditions of type
$\Theta_j$, for $j=1,\ldots,t$, (provided the codimensions of all the imposed
conditions add up to $\dim \M_{0,n}(X,\beta)$), where $n=\sum a_i + \sum b_j$).
Assume the tangency conditions satisfy assumption~\ref{assumptionWV}.

\bigskip

Since the genus zero spaces $\M_{0,\emptyset}(X,\beta)$ satisfy
assumption~\ref{assumptionM}, and the general map is irreducible, we have:
\begin{prop}\label{N=til}
	$
	N^{(0)}_\beta\!(a_{p+1},\ldots,a_r,b_1,\ldots,b_t)
	=
	\til N^{(0)}_\beta\!(a_{p+1},\ldots,a_r,b_1,\ldots,b_t).
	$
	\qed
\end{prop}

The corresponding potential is related to the tangency quantum potential by 
the linear coordinate change of Lemma~\ref{GGamma}, so an easy application of
the chain rule translates Theorem~\ref{PDE} into
\begin{satz}%
\textsc{Topological recursion for genus $0$ characteristic numbers.}%
\label{PDE-char}
	The following differential equations hold for the generating function of the
	genus zero characteristic numbers.
\begin{eqnarray*}
	G^{(0)}_{v_k u_i u_j} &= & \sum_{q=0}^r \rho_{qk} \; G^{(0)}_{u_q u_i u_j}
  \\
	 &+&\sum_{q=0}^r \sigma_{qk} \bigg(
G^{(0)}_{u_q(u_i u_j)} - G^{(0)}_{(u_q u_i)u_j}
	- G^{(0)}_{(u_q u_j)u_i} +
	\sum_{e,f} G^{(0)}_{u_q u_e} \, \gamma^{ef}(\fat{v}) 
	\, G^{(0)}_{u_f u_i u_j}
	\bigg).
\end{eqnarray*}
Here $(\gamma^{ef}(\fat{v}))$ denotes the matrix $(\gamma^{ef})$ with $y_i$
substituted by $\sum_j \sigma_{ij} v_j$.  The coefficients $\rho_{ij}$ and 
$\sigma_{ij}$ are those describing the tangency conditions, cf.~(\ref{rhosigma}).
\qed
\end{satz}

\begin{eks}\textsc{Projective space.}\label{Pr}
	Consider $X=\P^r$ with $h\df c_1(\OO(1))$ and the natural cohomology basis
	$T_i = h^i$.  Consider the conditions $\Theta_j = \mtau_0(T_{j+1}) +
	\mtau_1(T_j)$ of being tangent to a hyperplane $H$ along a specified
	co\-dimension-$j$ plane contained in $H$ ($j=1,\ldots,r$), together with the
	conditions $\Omega_j = \mtau_0(T_j)$ of being incident to a codimen\-sion-$j$
	plane ($j=2,\ldots,r$).  Let $N_d(\fat{a},\fat{b})$ denote the number of
	rational curves of degree $d$ satisfying $\fat{a}$ such incident conditions
	and $\fat{b}$ such tangency conditions.  Then by Proposition~\ref{PDE-char},
	the corresponding generating
	function $G(\fat{u},\fat{v})$ (as defined in \ref{gen-function}) satisfies the 
	differential equations
$$
G_{v_k u_1 u_1} =  G_{u_1 u_1 u_{k+1}} +
G_{u_2 u_k} - 2 \, G_{u_{k+1} u_1} +
	\sum_{e,f}  G_{u_k u_e} \, \gamma^{ef}(\fat{v}) \, G_{u_f u_1 u_1} .
$$
Here the matrix $(\gamma^{ef}(v))$ is obtained from $(\gamma^{ij}) =
(\phi^{ij}(2\fat{y}))$ by substituting $v_i$ for $y_i$. 
The case $X=\P^2$ is treated in more depth in the next section.
\end{eks}

\begin{eks}\label{G24}\textsc{The Grassmannian $\Gr(2,4)$.}
	Let $X=\Gr(2,4)$, the Grassmannian of lines in $\P^3$.  Take the basis of 
	Schubert varieties, such that $T_0$ is the fundamental class, $T_1$ is the 
	hyperplane section (under the Pl\"ucker embedding $\Gr(2,4)\into\P^5$), $T_2$ 
	and $T_3$ are the two Schubert varieties of dimension 2, $T_4$ is dual to 
	$T_1$ (it's a line), and $T_5$ is the class of a point.	 
	Let $N_d(b,a_2,a_3,a_4,a_5)$ be the number of rational curves of class
	$d\cdot T_4$ that are tangent to $b$ subvarieties of class $T_1$ and incident
	to $a_i$ subvarieties of class $T_i$, for $i=2,3,4,5$.  The tangency
	condition is of class $\Theta = \mtau_0(T_2) + \mtau_0(T_3) +
	\mtau_1(T_1)$.  Let $G(v, u_1,\ldots,u_5)$ denote the corresponding
	generating function, as in \ref{gen-function}.  By
	Proposition~\ref{PDE-char}, the differential equation for $G$ is
$$
G_{v u_1} = G_{u_1 u_2} + G_{u_1 u_3} - G_{u_2} - G_{u_3}
          + \brok{1}{2}
			 \sum_{e,f} G_{u_1 u_e} \, \gamma^{ef}(v) \, G_{u_f u_1}
$$
which is the integrated form, cf.~\ref{integration}.
Here the matrix is
$$
(\gamma^{ef}(v)) = 
\begin{pmatrix}
0 &      0 &      0 &      0 &      0 &      1        \\
0 &     0 &     0 &     0 &     1 &     2 v        \\
0 &    0 &    1 &    0 &    2 v &    2 v^2        \\
0 &    0 &    0 &    1 &    2 v &    2 v^2        \\
0 & 1 & 2 v & 2 v & 4 v^2  & \brok{8}{3} v^3          \\
1 &2 v &2 v^2  &2 v^2  &\brok{8}{3} v^3  &\brok{4}{3} v^4      
\end{pmatrix} .
$$
\end{eks}

\section{Characteristic numbers of plane curves \\ (of genus $0$, $1$, $2$)}

In this section we specialize to the case $X=\P^2$.  In genus $0$, the virtual
characteristic numbers agree with enumerative geometry.  In genus $1$ and $2$,
simple correction formulas are described.  The outlook for general genus appears
quite complicated due to difficult multiple cover contributions.

\subsection{Notation and main result}

\begin{blanko}{The characteristic number potential.}
	Let $N^{(g)}_d\!(a,b,c)$ denote the number of irreducible plane curves of genus
	$g$ and degree $d$ which are incident to $a$ general points, are tangent to
	$b$ general lines, and are tangent to $c$ general lines at a specified point
	on the line.  (This last condition will be referred to as a {\em flag}
	condition.)  Define the number to be zero if not $a+b+2c = 3d+g-1$.

Put
$$
G^{(g)}(s,u,v,w) = \sum_{d>0} \exp(ds) 
\sum_{a,b,c} \frac{u^a}{a!}\frac{v^b}{b!}\frac{w^c}{c!} \; N^{(g)}_d\!(a,b,c) .
$$
\end{blanko}

\begin{blanko}{Virtual characteristic numbers and their potentials.}
	As in \ref{virtualN}, define virtual characteristic numbers
\begin{equation}\label{virdeff}
	\tilde{N}^{(g)}_d(a,b,c) \df \Brac{ \;
	\big( \mtau_0(h^2) \big)^a \;
	\big( \mtau_0(h^2) \! + \! \mtau_1(h) \big)^b \;
	\big( \mtau_1(h^2) \big)^c \;
	}\!_{g,d}^{\P^2} ,
\end{equation}
and form the corresponding generating function
$$
\til G^{(g)}(s,u,v,w) = \sum_{d>0} \exp(ds) \sum_{a,b,c} 
\frac{u^a}{a!}\frac{v^b}{b!}\frac{w^c}{c!} \; \til N^{(g)}_d\!(a,b,c) .
$$
It is related to the tangency quantum potential $\Gamma^{(g)}$ by 
\begin{align*} \phantom{mmmmmmm}
\til G^{(0)}(s,u,v,w) &= \Gamma^{(0)}(x_1,x_2,y_1,y_2)  ,\\
\til G^{(1)}(s,u,v,w) &= \Gamma^{(1)}(x_1,x_2,y_1,y_2) 
+ \frac{1}{8} x_1  , \\
\til G^{(2)}(s,u,v,w) &= \Gamma^{(2)}(x_1,x_2,y_1,y_2) 
+ \frac{1}{960} y_1  , 
\quad \text{ and} \\
\til G^{(g)}(s,u,v,w)& = \Gamma^{(g)}(x_1,x_2,y_1,y_2)  
\phantom{+ \frac{1}{960} y_1}
\qquad \text{ for }g\geq 3 ,
\end{align*}
subject to the change of variables:
\begin{equation}
\label{cvar}
x_1=s , \qquad  x_2 = u + v ,  \qquad y_1 = v , \qquad y_2 = w .
\end{equation}
%
\end{blanko}

\begin{blanko}{}
Substituting this into the deformed metric (cf.~Example~\ref{gammaP2}), yields
the matrix
$$
\begin{pmatrix}
	\phantom{0}0\phantom{0} & \phantom{0}0\phantom{0}& \phantom{0}1\phantom{0} \\
	0 & 1 & 2v \\
	1 & 2v & 2v^2 + 2w
	\end{pmatrix} 
	$$
which comes up naturally together with partial derivatives whenever reducible
curves are in play.  We define two differential operators corresponding to the
last two rows of the matrix, (the ``line operator'' and the ``point operator''):
\begin{align}
L &\;\df\; \phantom{2v}\frac{\d}{\d s} + 2v \frac{\d}{\d u}\\
P &\;\df\; 2v\frac{\d}{\d s} + (2v^2+2w) \frac{\d}{\d u} .
\end{align}	
\end{blanko}


\begin{satz} 
\label{mainrel}
The virtual function $\tilde{G}^{(g)}$ is
related to enumerative geometry in genus $0$, $1$, and $2$ by the
following formulas.
\begin{enumerate}
\item[0.] $\tilde{G}^{(0)}= G^{(0)}$.
\item[1.] $\tilde{G}^{(1)}=G^{(1)}-\frac{1}{24}P G^{(0)}+
E$.
\item[2.] $\tilde{G}^{(2)}= G^{(2)} - \frac{1}{24}P G^{(1)} 
+\frac{1}{2}(\frac{1}{24}P)^{2}G^{(0)}
+ \big( H_s \cdot L G^{(0)} + H_u \cdot P G^{(0)} \big) +
Q_2 + Q_3$.
\end{enumerate}
\end{satz}
Here $E = \brok{1}{2} e^{2s}\big(
\frac{1}{2}\frac{v^6}{2!2!2!} + \frac{2u}{1!}\frac{v^5}{2!3!}
+ \frac{(2u)^2}{2!}\frac{v^4}{4!} + \frac{v^4}{2!2!}\frac{w}{1!}
+ \frac{2u}{1!}\frac{v^3}{3!}\frac{w}{1!} +
\frac{v^2}{2!}\frac{w^2}{2!}
\big)$ is the generating function for the elliptic double covers of a line
in $\P^2$, and
$H = \brok{1}{2} e^{2s}\big(
\frac{1}{2}\frac{v^8}{2!2!4!} + \frac{2u}{1!}\frac{v^7}{2!5!}
+ \frac{(2u)^2}{2!}\frac{v^6}{6!} + \frac{v^6}{2!4!}\frac{w}{1!}
+ \frac{2u}{1!}\frac{v^5}{5!}\frac{w}{1!} +
\frac{v^4}{4!}\frac{w^2}{2!}
\big)$ is the generating function for genus $2$ double 
covers of a line in $\P^2$.
(To be defined properly in the proof).
The terms $Q_2$ and $Q_3$ 
correspond to degenerate genus $2$ curves (of degree 2 or 3).
We will not explicitly compute these terms, since anyway we know there are
no immersions in these low degrees.

\bigskip


Theorem~\ref{mainrel} is proven in the next three subsections.  The strategy of
proof is straightforward: first, in~\ref{spurious} a geometric argument shows
the intersection locus corresponding to a virtual characteristic integral
(\ref{virdeff}) is a union of the enumerative solutions and excess loci.  The
arguments of Proposition~\ref{enum} show that the desired solutions are
isolated and count with multiplicity one.  Second, in order to conclude our
result, we just need to identify the contributions from the other loci, which we
can do explicitly in genus $1$ and $2$ by an argument involving the virtual
fundamental class.  These arguments are given in \ref{genus1} and \ref{genus2}.

\begin{blanko}{Genus 0.}
	The genus $0$ part of Theorem~\ref{mainrel} was proved in \ref{N=til}.
	Let us take the opportunity here to spell out the topological recursion 
	relations.
	Theorem~\ref{PDE-char} gives two differential equation for $G^{(0)}$ 
	(cf.~Example~\ref{Pr}).
	The first one reads, in its integrated form
	(cf.~\ref{integration}),
\begin{equation}\label{eqGP2}
G^{(0)}_{v s} =  G^{(0)}_{u s} - G^{(0)}_{u}
+\brok{1}{2} \big(  G^{(0)}_{s s}\cdot LG^{(0)}_{s}
+ G^{(0)}_{u s}\cdot P G^{(0)}_{s} \big).
\end{equation}
This equation determines the simple characteristic numbers $N^{(0)}_d\!(a,b,0)$
from the Gromov-Witten invariants.  It was first found
by the third named author in 1997 (cf.~\cite{Pand:mail}), and later, with
different methods, by R.~Vakil~\cite{Vakil:9803}.

Together with the second equation,
\begin{equation*}
G^{(0)}_{w s s} = G^{(0)}_{u u}
+ \big( G^{(0)}_{u s} \cdot L G^{(0)}_{ss}
+ G^{(0)}_{u u} \cdot P G^{(0)}_{ss} \big) .
\end{equation*}
it provides a transparent and effective way to compute the characteristic
numbers $N^{(0)}_d\!(a,b,c)$ from the Gromov-Witten invariants.
\end{blanko}

\subsection{Identifying spurious components}\label{spurious}

Our goal is to be able to give an exhaustive list of irreducible components of
the map space which contribute to our virtual characteristic number.  Since the
presence of marks multiplies the number of components, we perform the
dimension counts only in the markless situation.  

A necessary condition for a component of the map space to contribute is that it
have dimension greater than or equal to the expected dimension.  Such components
arise in connection with multiple cover maps or contracted tails.  The
condition is not sufficient however, because the maps in such a component may
not be able to satisfy all the conditions we impose (according to the expected
dimension).  In the terminology of Vakil~\cite{Vakil:9709003}, we need only care
about components of the map space whose {\em intersection dimension} is equal to
or greater than the expected dimension of the whole space.  The intersection
dimension of a component of the map space is defined as the maximum number $i$
such that there exist $a + b + 2c = i$ and $\incidenceclass^a \cp
\tangencyclass^b\cp\flagclass^c \neq 0$ in the operational Chow ring of the
component.  (Here $\flagclass$ stands for the class of maps tangent to a given
line at a specified point.)  Obviously the intersection dimension is bounded
above by the actual dimension.

\begin{blanko}{Geography of the map space.}\label{geography}
	To get a handle on the collection of irreducible components of the space of
	maps, we introduce the {\em label} of a stable map.  Given a stable map
	$\mu:C \rarr \P^2$ we define its label to be the dual graph of the curve $C$
	together with the quadruple $[g_i,h_i,e_i,k_i]$ associated to each vertex;
	here $g_i$ is the genus of the corresponding curve $C_i$; $h_i$ is the
	geometric genus of the image of $C_i$; $e_i$ is the degree of the image of
	$C_i$; and $k_i$ is the degree of $\mu$ restricted to $C_i$ onto its image.
	(If $e_i=0$ then we define $k_i$ and $h_i$ to be zero too.)
	We define the label of a component of the map space to be the label of the
	map parametrized by a general point of the component.  Conversely, we can
	associate to each label the locus in the map space parameterizing maps with
	that label.  Thus the labels stratify the moduli space.
\end{blanko}

We will now compute the dimension of the family of given label.

\begin{blanko}{Irreducible maps.}
	It is straightforward to compute the dimension of the family of a given label
	$[g,h,e,k]$.  The map is then a composition $C \pil C' \pil \P^2$, where the
	first map is a $k$-sheeted covering of a curve of genus $h$ by a curve of
	genus $g$, and the second map is an immersion.  The dimension of this family
	is easily computed as the sum of the dimension of the Hurwitz scheme and the
	Severi variety, so it is
\begin{equation}\label{irrdimformula}
	(2g-2) - k(2h-2) \ \ + \ \ 3e + h - 1  .
\end{equation}
Observe that this is also the intersection dimension, since for example the
image curve can satisfy $3e+h-1$ incidence conditions, and then the $(2g-2) -
k(2h-2)$ ramification points can account for a tangency condition each.

Requiring this number to be greater than or equal to the expected dimension $3ke
+ g - 1$, can be written as the inequality
\begin{equation}\label{inequality}
g-h + (1-k)(3e+2h-2) \geq 0 .
\end{equation}
If strict equality holds, the left hand side is the excess dimension.  It is
easy to see that (except for the trivial solution $h=g$, $k=1$ corresponding to
immersions) it is necessary to have $h<g$, and that for fixed $g$ and $h$,
there is only a finite set of solutions (which are small values of $e$ and $k$).
%
%
%
%
%
\end{blanko}
\begin{blanko}{Reducible maps --- no contraction.}\label{reducible maps}
	When there are no contracted components of the map, the dimension of the
	family of given label does not depend of the structure of the graph, but only
	on the set of quadruples $g_{i}, h_{i},e_{i},k_{i}$, associated to the
	vertices.  The dimension of the family is just the sum of the dimensions of
	the space parameterizing the maps restricted to each component.
	Observe that we will always be able to assemble the collection of curves into
	a stable map, since there will always be intersection points of the various
	image curves which we can use to glue the curves together.  Moreover,
	generically there will be only finitely many such points, so this choice does
	not affect the dimension of the component of the map space.

	The important thing to note is that the expected dimension formula has a
	super\-additivity property.  If $C=C'\cup C''$ is a union of curves of
	degrees $d'$ and $d''$ and genera $g'$ and $g''$, then $C$ has degree
	$d'+d''$ and genus $g'+g''$.  So its expected dimension is
\[
3(d'+d'') + (g'+g'')-1=(3d'+g'-1)+(3d''+g''-1)+1.
\]
	From this formula it follows that in order for a reducible curve to move in
	dimension greater than or equal to the expected dimension, one of its
	components must move in dimension strictly higher than expected.
\end{blanko}

\begin{blanko}{Contracted tails.}\label{contractedtails}
	A contracted component can contribute large amounts to the dimension of a
	component of the map space, due to the freedom of varying the moduli of the
	contracted curve.  However, this freedom is irrelevant for the sake of
	satisfying the conditions we are imposing.
	
	The only way a contracted component can satisfy a condition of being tangent
	to a given line is when it maps to a point on that line.  If the component is
	attached to the rest of the curve in only one point (in which case we call it
	a tail), then this possibility represents one degree of freedom, namely:
	after the honest part of the curve has satisfied all the conditions it can
	(according to the count of \ref{reducible maps}), one more tangency condition
	can be satisfied, by choosing the attachment point in such a way that it maps
	to the given line.  (Alternatively, we could begin by requiring the
	contracted tail to map to the intersection of two given lines or to the
	point of a flag condition, thus satisfying two extra conditions, but this
	imposes an extra incidence condition on the rest of the curve --- since in
	order to glue, it needs to pass through the point.)  Either way, we see that
	a contracted tail contributes exactly one to the intersection dimension of
	the component of the map space, compared to the intersection dimension
	computed in the previous paragraphs ignoring this contracted tail.
	
	On the other hand, if the contracted component is attached to two or more
	points then its image point is a node of the image curve, and since the image
	curve has only a finite number of nodes there is no freedom left to satisfy
	further conditions.  So in this case there is no contribution to the
	intersection dimension.  In particular, by stability, we see that a
	contracted curve of genus $0$ never contributes to the intersection
	dimension.
\end{blanko}
	
In summary, we get (an upper bound on) the intersection dimension of a component
of the map space in terms of its label, namely the sum of the dimensions given by
formula~(\ref{irrdimformula}) for each noncontracted component, plus the total
number of contracted tails.
	
Now we are in a position to find all labels in low genus whose associated
intersection dimension is greater than or equal to the expected dimension.

\subsection{Genus 1}\label{genus1}

\begin{prop}\label{genus1geo}  The following is a complete list
of relevant combinatorial
types of genus $1$ curves: 
\begin{punkt-i}
\item [(i)] For all $d\geq 3$, an immersed degree $d$ elliptic curve.
\item [(ii)] For $d=2$, an elliptic double cover of a line.
\item [(iii)] For all $d\geq 1$, a rational degree $d$ curve with a contracted 
elliptic tail.
\end{punkt-i}
\end{prop}

\begin{dem}
	It is clear that these types occur, we need to rule out the existence of
	other types.  First observe that no irreducible maps move in dimension higher
	than expected and that those moving in expected dimension are the immersions
	and the double covers of lines.  Indeed, to get solutions to the
	inequality~(\ref{inequality}), we need $g=1$ and $h=0$,
	whereafter it reads
	$$
	1 + (1-k)(3e-2) \geq 0.
	$$
	Clearly the only solutions are $k=1$, $e=e$ (type (i))
	and $k=2$, $e=1$ (type (ii)).
	
	If there is a contracted tail, it must be of genus $1$ by the observation of
	\ref{contractedtails}.  Then the rest of the curve must have genus $0$ and
	must move in the expected dimension; hence it is an immersion.  This is type
	(iii).
\end{dem} 

In order to prove the genus $1$ case of \ref{mainrel} we need to determine the
contributions of curves of types (ii) and (iii) to the descendant generating
series.

\begin{blanko}
{The type (ii) contribution} is $E = \brok{1}{2} e^{2s}\big(
\frac{1}{2}\frac{v^6}{2!2!2!} + \frac{2u}{1!}\frac{v^5}{2!3!} +
\frac{(2u)^2}{2!}\frac{v^4}{4!} + \frac{v^4}{2!2!}\frac{w}{1!} +
\frac{2u}{1!}\frac{v^3}{3!}\frac{w}{1!} + \frac{v^2}{2!}\frac{w^2}{2!} \big)$. 
This presents no difficulty, since the moduli space of these curves has the
expected dimension.  We just need to count the solutions.

It is impossible for a line to meet more than 2 points, so we get 
contributions only for the 6 characteristic numbers 
$\til N^{(1)}_2\!(0,6,0)$,
$\til N^{(1)}_2\!(1,5,0)$, 
$\til N^{(1)}_2\!(2,4,0)$,
$\til N^{(1)}_2\!(0,4,1)$,
$\til N^{(1)}_2\!(1,3,1)$, and
$\til N^{(1)}_2\!(0,2,2)$ 
which correspond exactly to the six terms in the polynomial $E$.

Let us count those with six tangency conditions: we need to choose two pairs of
the six given lines (there are $\frac{1}{2}\binom{6}{2,2}=45$ ways to do that),
and then draw the solution curve as the unique line joining the two
corresponding intersection points.  Then the ramification points of the map are
completely determined: two of them must be the inverse images of the two points,
and the other two must be the intersection points of the image line with the
remaining two given lines.  This explains the term
$\frac{1}{2}\frac{v^6}{2!2!2!}=45\frac{v^6}{6!}$.  

(Throughout the polynomial, the variable $u$ appears together with a factor $2$,
since for each incidence condition, there are two ways of putting the
corresponding mark.  All solutions would be counted twice this way (or otherwise
possess an automorphism of order $2$); this is corrected for by the coefficient
$\frac{1}{2}$ in front of the whole polynomial.)
%
\end{blanko}

\begin{blanko}{Contribution from type (iii).}\label{virtualargument}
	The more interesting contributions arise from the curves with a contracted
	elliptic tail.  There are three types of solutions here.  Either the
	contracted tail can account for one tangency condition (by mapping to a
	point on the given line), or it can satisfy two such tangency conditions (by
	mapping to the intersection of the two lines), or it can satisfy a flag
	condition (by mapping to the point of the flag).  The analysis is the same in
	each case so we give only the first case.
	
The component of genus zero satisfies all of the conditions except for one
tangency condition, and the genus $1$ tail is contracted over one of the $d$
intersection points of the genus zero curve with the remaining line.  To
compute the contribution of these solutions, we first have the combinatorial
factor of $bdN^{(0)}_d\!(a,b-1,c)$ accounting for choosing which of the $b$
lines the contracted curve maps to, which of the $N^{(0)}_d\!(a,b-1,c)$
rational curves through the other points we choose, and which of the $d$
points of intersection of that curve with the given line will be the
point of attachment for the contracted elliptic tail.  However, there remains
an entire $\M_{1,1}$ of solutions corresponding to the choice of $j$-invariant
of the elliptic curve over which we must integrate the virtual class.

We will compute this as an integral on the one-pointed space of maps having one
contracted elliptic tail and one rational component of degree $d$.  Note that
this space has two components: one where the marked point is on the rational
curve and one where it is on the elliptic curve.  We compute the contributions
from these two components separately.  On the first component of the space the
cycles intersect transversally in a substack isomorphic to $\M_{1,\{x\}}$ ($x$
is the gluing mark on the elliptic curve).  Thus, the contribution is the degree
of the virtual class restricted to this locus.  Since the unpointed space is
smooth near this locus, it is simply the first Chern class of the obstruction
bundle $\EE$, which fits into the following exact sequence.
\[
0 \rarr \LL\dual_x \rarr \Hodge\dual \sumo \Hodge\dual \rarr \EE \rarr 0.
\]
Here $\Hodge$ is the Hodge bundle.  Since $\Hodge$ is isomorphic to $\LL_x $ on
the one-pointed space $\M _{1,\{x\}}$, we conclude that the restriction of the
virtual class to this locus is $-\mpsiclass_x$.  Hence the contribution of this
component is $-1/24$ (cf.~(\ref{specialdilaton})).\nocite{Graber:9808}

On the other component of the space, we see that in addition to the virtual
class, we need to account for an excess intersection.  This is straightforward. 
Once we impose the condition that the marked point meet a line, the entire
contracted tail is already forced to lie over the line.  Hence, the second
condition $\evalmap_1\upperstar (h) + \mpsiclass_1$ is the excess class.  As it
is impossible for this curve to meet another general line, the contribution from
$\evalmap_1\upperstar (h)$ is zero, and we are left with $\mpsiclass_1$ which we
must integrate against the virtual class.  The solution locus is isomorphic to
$\M _{1,\{p_1,x\}}$; by Lemma~\ref{restr-mpsi}, $\mpsiclass_1$ restricts to give
$\mpsiclass_1 + \deltaclass_{x1}$.  The virtual class is again $c_1(\Hodge\dual)
= -\mpsiclass_x$, so in the end we find
$$
\int_{\M_{1,\{p_1,x\}}} \!\!\!\!\!\!(\mpsiclass_1 + \deltaclass_{x1}) \cp
(-\mpsiclass_x) \; = \; - \int_{\M_{1,\{p_1,x\}}}
\!\!\!\!\!\!\mpsiclass_1\mpsiclass_x - \int_{\M_{1,\{x\}}}
\!\!\!\!\!\!\mpsiclass_x \; = \; -\frac{1}{24},
$$
by Lemma~\ref{push-delta}, the dilaton equation~(\ref{dilaton}), and 
(\ref{specialdilaton}) again.

In total, then, we see that the contribution of these curves to the virtual
characteristic number is $-\frac{1}{24}\, 2bd \, N^{(0)}_d\!(a,b-1,c).$
Similarly, we find that the locus when the contracted elliptic tail lies over
the intersection of two of the lines gives a contribution of $-\frac{1}{24} \, 4
\binom{b}{2} \, N^{(0)}_d\!(a+1,b-2,c)$, and finally, in the case where the
contracted elliptic tail accounts for a flag condition, the contribution is
$-\frac{1}{24} \, 2c \, N^{(0)}_d\!(a+1,b,c-1)$.  Combining these contributions
gives us exactly the formula in the genus $1$ case of Theorem~\ref{mainrel}. 
\qed
\end{blanko}

\begin{blanko}{Topological recursion relations for the genus $1$ characteristic
numbers.}
	In order to derive a TRR for the characteristic numbers, it is convenient to 
	include the degree $2$ term $E$, so we let now $G^{(1)}$ denote what was called
	$G^{(1)}+E$ in Theorem~\ref{mainrel}.
%
%
	Now the chain rule applied to \ref{TRR1} gives differential equations for
	$\til G^{(1)}$, which combined with
%
Theorem~\ref{mainrel} (and the genus $0$ TRR backwards
twice) yields the following TRR for the genus $1$ characteristic number 
potential:
\begin{multline*}
G^{(1)}_v = G^{(1)}_u 
+ \big( 
L G^{(0)}_{s} \cdot  G^{(1)}_s 
+ P G^{(0)}_{s} \cdot  G^{(1)}_u 
\big)
 \\ +
\frac{1}{24}\bigg( LG^{(0)}_{ss} + PG^{(0)}_{us} - 2 L G^{(0)}_s
+ 2 G^{(0)}_s - G^{(0)}_{ss} \
 - 2vG^{(0)}_{vs} - (2v^2 + 2w) G^{(0)}_{ws}\bigg) .
\end{multline*}
This equation has also been established by R.~Vakil~\cite{Vakil:9803} by
degeneration methods and formulas from \cite{Pand:9504}.  His equation however
includes potentials corresponding to characteristic numbers of maps with a
double point on a given line or at a given point.

Similarly, there is an equation for taking away a flag condition:
\begin{multline*}
G^{(1)}_w = \big( 
L G^{(0)}_{u} \cdot  G^{(1)}_s 
+ P G^{(0)}_{u} \cdot  G^{(1)}_u 
\big)
 \\ +
\frac{1}{24}\bigg( LG^{(0)}_{us} + PG^{(0)}_{uu} - 2 L G^{(0)}_u
+ 2 G^{(0)}_u - G^{(0)}_{us} \
 - 2vG^{(0)}_{vu} - (2v^2 + 2w) G^{(0)}_{wu}\bigg) .
\end{multline*}
Together these two equations determine all the numbers $N^{(1)}_d\!(a,b,c)$ from
the Gromov-Witten invariants $N^{(1)}_d = N^{(1)}_d\!(3d,0,0)$ and the genus $0$
characteristic numbers $N^{(0)}_d\!(a,b,c)$.  The Gromov-Witten invariants are
determined by the recursion of Eguchi-Hori-Xiong~\cite{Eguchi-Hori-Xiong},
proved in Pandharipande~\cite{Pand:9705} using the relation of
Getzler~\cite{Getzler:9612}.
\end{blanko}

\subsection{Genus 2}\label{genus2}

The proof of Theorem~\ref{mainrel} proceeds in essentially the same way for
genus $2$.  First we identify the relevant components of the map space.

\begin{prop}\label{genus2geo}
The following is a complete list of relevant combinatorial types
of genus $2$ maps with $d\geq 4$:
\begin{punkt-i}
\item  Immersed genus $2$ curve.
\item [(ii)] Immersed genus $1$ curve with a 
contracted elliptic tail.
\item [(iii)] Immersed rational curve with two contracted 
elliptic tails.
\item [(iv)] Immersed degree $(d-2)$ rational curve attached to a genus $2$ double
cover of a line.
\end{punkt-i}
\end{prop}
(In lower degree, there are further three special cases.  In degree $2$: a
genus $2$ double cover of a line (they move in excessive dimension), and an
elliptic double cover of a line with a contracted elliptic tail (expected
dimension).  And in degree $3$: a genus $2$ triple cover of a line (they move in
expected dimension).  


\begin{dem}
The proof is essentially the same as the genus $1$ case, so we just
indicate the new points.  First, the possibility of a genus $2$ contracted
tail is ruled out by the fact that the rest of the curve would
then have to be genus $0$ which could only move in dimension $3d-1$.  The
contracted tail contributes only 1 to the intersection dimension,
which is not enough to give the expected dimension of $3d+1$.
The only other difference is that the genus $2$ double covers of lines move
in dimension 1 greater than expected.  Hence, they can occur as a component
of a reducible curve moving in the expected dimension, accounting for
the 4th type.  No other irreducible curve moves in greater than expected
dimension, so these are the only interesting reducible curves.
\end{dem} 

\begin{blanko}{Contribution from type (ii) and (iii).}	
	If we have just one contracted elliptic tail, then all of the analysis we did
	in the genus $1$ case still goes through.  The solution loci are still either
	$\M_{1,1}$ or $\M_{1,2}$ and the obstruction theory is the same.  We conclude
	that we get contributions of $ -\frac{1}{24}\big(  2bd \, N^{(1)}_d\!(a,b-1,c)
	+ 4 \binom{b}{2}\, N^{(0)}_d\!(a+1,b-2,c) + 2  c \,
	N^{(0)}_d\!(a+1,b,c-1) \big)$, 
	which in terms of the potentials amounts to $-\frac{1}{24}P G^{(1)}$.
	
	If we have two contracted tails, then each of them can be used in either of
	the three ways described (satisfying one or two tangency condition, or one
	flag condition), so there are six combinatorial types of contributions here. 
	The locus of solutions is either $\M_{1,1} \times \M_{1,1}$, $\M_{1,1} \times
	\M_{1,2}$, or $\M_{1,2} \times \M_{1,2}$.  In each case, the obstruction
	theory is the product obstruction theory.  This is intuitively clear, since
	the deformations of the two contracted tails are unrelated.  (It is also easy
	to verify directly that the arising excess bundles and obstruction bundles
	naturally decompose into products of the ones that occur in the
	one-contracted-tail case.)
	The contribution from this type amounts to
	$$
	\frac{(\frac{1}{24}P)^2}{2!} \;G^{(0)},
	$$
	as claimed.  (One can verify this by writing out the six terms,
	\begin{equation*}\begin{split}
	\frac{1}{24^2}\bigg(4d\smallbinom{b}{2}N^{(0)}_d\!(a,b-2,c)
	+ 8d\smallbinom{b}{2,1}N^{(0)}_d\!(a+1,b-3,c)
	+ 8d\smallbinom{b}{2,2}N^{(0)}_d\!(a+2,b-4,c) \\
	+ 4bcN^{(0)}_d\!(a+1,b-1,c-1)
	+ 8c\smallbinom{b}{2}N^{(0)}_d\!(a+2,b-2,c-1)
	+ 4\smallbinom{c}{2}N^{(0)}_d\!(a+2,b,c-2)
	\bigg) ,
	\end{split}\end{equation*}%
and translating back into the language of potentials).
%
\end{blanko}

\begin{blanko}{Contribution from type (iv).}\label{iv}
	The only remaining point then, is to compute the contributions coming 
from the
curves of type (iv).  These move in the expected dimension, so we need only 
count them.  This is reasonably straightforward, but there are many
types of solutions.

The genus $2$ double covers of a line can satisfy $8$ conditions, of which at 
most two can be incidence conditions.  The generating function for the numbers 
of genus $2$ double covers is 
$$
H = \brok{1}{2} e^{2s}\big(
\frac{1}{2}\frac{v^8}{2!2!4!} + \frac{2u}{1!}\frac{v^7}{2!5!}
+ \frac{(2u)^2}{2!}\frac{v^6}{6!} + \frac{v^6}{2!4!}\frac{w}{1!}
+ \frac{2u}{1!}\frac{v^5}{5!}\frac{w}{1!} +
\frac{v^4}{4!}\frac{w^2}{2!} 
\big) .
$$
This follows from the same arguments as those giving the potential $E$ of
elliptic double covers.  It is important to note that the coefficients are the
actual number of such maps, but that they are not Gromov-Witten invariants! 
Indeed, the expected dimension of $\M_{2,0}(\P^2,2)$ is $7$, so all the
corresponding Gromov-Witten invariants are zero for dimension reasons (the
enumerative numbers are obtained integrating against the topological fundamental
class instead of the virtual one\ldots)

\bigskip

Now the solutions fall in two groups.  Either the genus $2$ component $C'$
satisfies eight conditions honestly or it satisfies seven.  If $C'$ satisfies
eight conditions honestly, then the rational curve $C''$ satisfies the remaining
conditions, either honestly or by attaching itself to $C'$ at the intersection
point of $\mu(C')$ with one of the given lines, thus accounting for tangency to
this line (this constitutes an extra incidence condition on $C''$).  The
contribution for these two cases is $H_s \cdot \big( G^{(0)}_s + 2v
G^{(0)}_u\big)$.  Here the subscript $s$ on $H$ corresponds to a factor $2$ in
front of the characteristic numbers, coming from the choice between the two
points of $C'$ that are possible for gluing.  Similarly, there is the choice
among $d-2$ points as attachment point on $C''$; this explains the subscript $s$
on $G^{(0)}$.  The coefficient $2$ reflects the fact that ``node on line''
counts twice as tangency: we are actually counting {\em marked} maps, and over
the locus of nodes there are two components, depending on which of the two
curves carries the mark.  The factor $v$ corresponds to a factor $b$ in front of
the characteristic number, accounting for the choice of which line we force the
node upon; and the subscript $u$ comes from the fact that the rational curve
acquires an extra incidence condition when forced to have the node mapping to
the intersection of the given line and $\mu(C')$.

If the genus $2$ curve $C'$ satisfies only seven conditions honestly, the
rational part must satisfy at least one extra tangency condition at the node
(which forces $C'$ to pass the point where this occurs thus imposing another
incidence condition on $C'$).  This case gives the contribution of $2v \!\cdot\! 
H_u \!\cdot\!  G^{(0)}_s$.  Alternatively, the rational part can satisfy two
tangency conditions by having the node on the intersection of two of the given
lines.  In this case both $C'$ and $C''$ acquire an extra incidence condition,
so that contributes $2v^2 \!\cdot\!  H_u \!\cdot\!  G^{(0)}_u$.  (Here the
multiplicity of a node accounting for two tangencies is $4$ (four ways to put
the two marks on the two components), and the choice of which two lines are used
gives a factor $\binom{b}{2}$ in front of the characteristic number and thus a
factor $\frac{v^2}{2!}$ in front of the potential.)  Finally, the node could
fall on the point of a flag condition, again imposing one extra incidence
condition on each of the components, thus giving a contribution of $2w \!\cdot\! 
H_u \!\cdot\!  G^{(0)}_u$.

This completes the proof of Theorem~\ref{mainrel}.
\end{blanko}

\section{Characteristic numbers of $\P^1\times\P^1$ (genus $0$ and $1$)}
\label{Sec:P1P1}
\setcounter{subsection}{1}
\setcounter{lemma}{0}

The techniques of the preceding section readily solve the characteristic number
problem for rational and elliptic curves in the quadric surface
$\P^1\times\P^1$.  In genus $0$, the result is immediate from
Proposition~\ref{N=til}.  In genus $1$, the correction terms involve the Hurwitz
numbers.


\begin{blanko}{Characteristic numbers of $\P^1$ --- the Hurwitz numbers.}
	For $X=\P^1$ (with $T_0$ = fundamental class, $T_1$ = class of a point), it
	is easy to prove that the invariants $N^{(g)}_d\!(b) = \brac{
	\mtau_1(T_1)^{b} }\!_{g,d}^{\P^1}$ are exactly the simple Hurwitz numbers
	(the number of $d$-sheeted genus $g$ coverings of the Riemann sphere simply
	ramified over $b=2d +2g-2$ given points).  Form the corresponding potential,
	\begin{equation}\label{Hurwitz}
		H^{(g)}(t,v) =
	\sum_{d>0} \exp(dt) \sum_b \frac{v^b}{b!} N^{(g)}_d\!(b) .
	\end{equation}
	The topological recursion relations \ref{PDE} and \ref{TRR1} directly give
	the well-known equations (see for example Vakil~\cite{Vakil:9803})
	\begin{align}
	H^{(0)}_{vt} &= v H^{(0)}_{tt} \!\cdot\! H^{(0)}_{tt} \\
	H^{(1)}_{v} &= 
2v H^{(0)}_{tt} \!\cdot\! H^{(1)}_t
+ \brok{1}{24} 2v \big( H^{(0)}_{ttt}
- H^{(0)}_{tt} \big).
\end{align}
\end{blanko}

%
%
%
%

\begin{blanko}{Set-up for $\P^1\times\P^1$.}
	Let $T_0$ be the fundamental class; let $T_3$ be the class of a point; and
	let $T_1$ and $T_2$ be the hyperplane classes pulled back from the two
	factors.  A curve of class $\beta$ is said to have bi-degree $(d_1,d_2)$,
	where $d_1=\int_\beta T_1$ and $d_2=\int_\beta T_2$.  A curve of bi-degree 
	$(1,0)$ is called a horizontal rule, and a curve of bi-degree $(0,1)$ a 
	vertical rule. Let
	$N^{(g)}_{(d_1,d_2)}(a,b,c)$ denote the characteristic numbers of irreducible
	curves in $\P^1\times\P^1$ of genus $g$ and bi-degree $(d_1,d_2)$ passing
	through $a$ general points, tangent to $b$ general curves of bi-degree
	$(1,1)$, and tangent to $c$ such curves at specified point.  Let
	$G^{(g)}(u_1,u_2,u,v,w)$ be the corresponding generating function ($u_1$ and
	$u_2$ being the formal variables corresponding to the partial degrees $d_1$
	and $d_2$).

The virtual classes corresponding to these three conditions are
\begin{align*}
\mtau_0(T_3) &\qquad \text{ (the incidence condition)}, \\
2\mtau_0(T_3) + \mtau_1(T_1) + \mtau_1(T_2) 
&\qquad \text{ (the tangency condition)},\\
\mtau_1(T_3) &\qquad \text{ (the flag condition)}.
\end{align*}
Let the virtual potential $\til G^{(g)}$ be defined as in \ref{gen-function}. 
Then we have $\til G^{g)}(u_1,u_2,u,v,w) =
\Gamma^{(g)}_+(x_1,x_2,x_3,y_1,y_2,y_3)$, with $x_1 = u_1$, $x_2=u_2$, $x_3= u +
2 v$; $y_1=v$, $y_2=v$, $y_3=w$.  For convenience, put also $$s\df u_1 + u_2,$$
the formal variable corresponding to $T_1 + T_2$.  Plugging these substitutions
into the matrix of the deformed metric yields
$$
\begin{pmatrix}
0 &      0 &      0 &            1   \\
0 &     0  &          1 &     2 v    \\
0 &      1 &    0 &    2 v        \\
1 &2 v&2 v    & 4 v^2 + 2w      
\end{pmatrix} .
$$
Define three differential operators corresponding to the three last lines of 
this matrix
\begin{align*}
L_1 &\;\df\; \phantom{2v}\brok{\d}{\d u_2} + 2v \brok{\d}{\d u}\\
L_2 &\;\df\; \phantom{2v}\brok{\d}{\d u_1} + 2v \brok{\d}{\d u}\\
P &\;\df\; 2v\brok{\d}{\d u_1} + 2v\brok{\d}{\d u_2} + (4v^2+2w) \brok{\d}{\d u} .
\end{align*}
\end{blanko}

\begin{blanko}{Genus $0$.}
In genus $0$, by Proposition~\ref{N=til}, the numbers encoded in the virtual
potentials are exactly the enumerative ones.
Proposition~\ref{PDE-char} gives 
\begin{align}
G_{vs} &= 2 G_{us} - 2 G_u + \frac{1}{2}\big( G_{s u_1} \!\cdot\!  L_1 G_{s}
+G_{s u_2} \!\cdot\!  L_2 G_{s} + G_{us} \!\cdot\! P G_{s} \big) \\
G_{wss} &=
2 G_{uu} + \big( G_{u u_1} \!\cdot\! L_1 G_{ss} +G_{u u_2} \!\cdot\! L_2 G_{ss}
+ G_{uu} \!\cdot\! P G_{ss} \big).
\end{align}
Together, these two equations constitute a transparent algorithm for computing
the genus $0$ characteristic numbers of $\P^1\times\P^1$ from the Gromov-Witten
invariants.  (Note that derivative
with respect to $s$ corresponds to a factor $d_1+d_2$ in front of the
characteristic number.)

Remark: for bi-degrees $(i,0)$ and $(0,i)$, with $i\geq 2$, the counted curves
are not immersions.
\end{blanko}

\subsection{Genus one}

\begin{lemma}\label{P1P1g1}
	The following is a complete list of relevant combinatorial types of genus $1$
	curves of bi-degree $(d_1,d_2)$ with $d_1>0$ and $d_2>0$:
\begin{punkt-i}
\item [(i)] An immersed  elliptic curve of bi-degree $(d_1,d_2)$.


\item [(ii)] An immersed rational curve of bi-degree $(d_1,d_2)$ with a
contracted elliptic tail.  

\item [(iii)] For all 
$2\leq i \leq d_1$, an elliptic $i$-sheeted cover of a horizontal rule union 
a rational curve of bi-degree $(d_1-i,d_2)$.  

\item [(iv)] For all 
$2\leq j \leq d_2$, an elliptic $j$-sheeted cover of a vertical rule union 
a rational curve of bi-degree $(d_1,d_2-j)$.  
\end{punkt-i}
\end{lemma}

\begin{dem}
	This follows readily by adapting the arguments of \ref{spurious},
	and Proposition~\ref{genus1geo}.
\end{dem}

\begin{blanko}{Multiple covers of a rule.}\label{rulecovers}
	The generating function for the genus $1$ covers of a horizontal rule
	is
	\begin{equation}
	I^{(1)}(u_1,u,v,w) = uH^{(1)}_{u_1} + (v^2+w)H^{(1)}_v ,
	\end{equation}
	where $H^{(1)}=H^{(1)}(u_1,v)$ is the Hurwitz potential, cf.~(\ref{Hurwitz}). 
	Indeed, the supporting rule for an $i$-sheeted map is fixed by either one
	incidence condition, two tangency conditions, or one flag condition.
	Once the supporting rule is fixed, the Hurwitz potential encodes the number
	of possible coverings.
	For the incidence condition, there are $i$ choices for the mark; this explains
	the factor $uH^{(1)}_{u_1}$.  For the case of $2i+1$ tangency conditions,
	the rule must pass through one of the two intersection point of two of the
	given curves.  This gives $2\!\cdot\!\binom{2i+1}{2}$ choices for the
	supporting rule, explaining the term $2\!\cdot\!\frac{v^2}{2!} H^{(1)}_v$. 
	Finally, for one flag condition and $2i-1$ tangency conditions, the flag
	fixes the supporting rule and translates into an extra $v$ condition on the
	covering of that rule.
		
	For coverings of a vertical rule we similarly find the generating function
	\begin{equation}
		J^{(1)}(u_2,u,v,w) = uH^{(1)}_{u_2} + (v^2+w)H^{(1)}_v.
	\end{equation}
\end{blanko}

\begin{prop}\label{111}
The virtual function $\tilde{G}^{(1)}$ for $\P^1\times\P^1$ is
related to enumerative geometry by the formula
$$
\til G^{(1)} = G^{(1)}-\frac{1}{24}P G^{(0)} \ +
\ I^{(1)}_{u_1} \!\cdot\! L_1 G^{(0)} + I^{(1)}_u \!\cdot\! P G^{(0)}
\ + \ J^{(1)}_{u_2} \!\cdot\! L_2 G^{(0)} 
 + J^{(1)}_u \!\cdot\! P G^{(0)} ,
$$
where $I$ and $J$ are the potentials defined in \ref{rulecovers}.
\end{prop}

\begin{dem}
	The term $-\frac{1}{24}P G^{(0)}$ appears for the same reason as the
	corresponding term in Proposition~\ref{genus1geo} and the proof is also the
	same --- see \ref{virtualargument}.  (Note however that the symbol $P$ stands
	for different operators in \ref{mainrel} and \ref{111}.)  The quadratic
	correction terms correspond to the reducible curves of type (iii) and (iv) in
	Lemma~\ref{P1P1g1}.  The potential for these reducible maps is found as in
	\ref{iv}.
\end{dem}

\bigskip

\small

\noindent
\textsc{Dept.~of Mathematics, 
Harvard University, 
Cambridge MA 02138}\\
{\em E-mail address:} graber@math.harvard.edu

\medskip

\noindent
\textsc{Dept.~of Mathematics, 
Royal Institute of Technology,
100 44 Stockholm, Sweden}\\
{\em E-mail address:} kock@math.kth.se

\medskip

\noindent
\textsc{Dept.~of Mathematics,
California Institute of Technology,
Pasadena, CA 91125}\\
{\em E-mail address:} rahulp@caltech.edu

\end{document}